\documentclass{gtart_h}  

\def\ifplaintex{\expandafter\ifx\csname documentclass\endcsname\relax}

\ifplaintex 
\else
\fi

\expandafter\ifx\csname beginpicture\endcsname\relax
\expandafter\ifx\csname documentclass\endcsname\relax
\input pictex \else
\input prepictex \input pictex \input postpictex \fi\fi

\def\gt{{\mathsurround=0pt\it $\cal G\mskip-2mu$eometry \&\ 
$\cal T\!\!$opology}}        

\def\gtp{{\mathsurround=0pt\it $\cal G\mskip-2mu$eometry \&\ 
$\cal T\!\!$opology $\cal P\!$ublications}}  

\def\lognumber#1{\def\thelognumber{#1}}
\def\volumenumber#1{\def\thevolumenumber{#1}}
\def\papernumber#1{\def\thepapernumber{#1}}
\def\volumeyear#1{\def\thevolumeyear{#1}}

\def\pagenumbers#1#2{\def\startpage{#1}\def\finishpage{#2}}
\def\published#1{\def\publishdate{#1}}
\def\proposed#1{\def\theproposer{#1}}
\def\seconded#1{\def\theseconders{#1}}
\def\received#1{\def\receiveddate{#1}}
\def\revised#1{\def\reviseddate{#1}}
\def\accepted#1{\def\accepteddate{#1}}

\def\asciiaddress#1{\def\theasciiaddress{#1}}
\def\asciiemail#1{\def\theasciiemail{#1}}

\long\def\asciiabstract#1{\long\def\theasciiabstract{#1}}

\let\\\par\let\thelognumber\relax
\let\thevolumenumber\relax\let\thepapernumber\relax
\let\thevolumeyear\relax\let\thesamplenumber\relax\let\startpage\relax
\let\finishpage\relax\let\publishdate\relax\let\receiveddate\relax
\let\reviseddate\relax\let\accepteddate\relax\let\theasciititle\relax
\let\theasciiauthors\relax\let\theasciiaddress\relax
\let\theasciiabstract\relax
\let\theasciiemail\relax\let\theshortauthors\relax\let\theshorttitle\relax

\long\def\maketitlep{   

\count0=\startpage

\gt\hfill      
\beginpicture
\setcoordinatesystem units <0.33truein, 0.33truein> point at 2.2 0.9
\setplotsymbol ({$\cal G$})
\plotsymbolspacing=9truept
\circulararc 315 degrees from 0 1 center at 0 0
\setplotsymbol ({$\cal T$})
\circulararc 315 degrees from 1 -1 center at 1 0
\endpicture
\break
{\small\ifx\thesamplenumber\relax 
Volume \else Sample
\fi\thevolumenumber\ (\thevolumeyear)
\startpage--\finishpage\nl
Published: \publishdate}
\vglue 0.5truein plus 0.4fil minus 0.1truein

{\parskip=0pt\leftskip 0pt plus 1fil\def\\{\par\smallskip}{\ifplaintex\large
\else\Large\fi\bf\thetitle}\par\medskip}   

\vglue 0pt plus 0.1fil 

{\parskip=0pt\leftskip 0pt plus 1fil\def\\{\par}{\sc\theauthors}
\par\medskip}

\vglue 0pt plus 0.1fil 

{\small\parskip=0pt\let\newline\\
{\leftskip 0pt plus 1fil\def\\{\par}{\sl\theaddress}\par}
\expandafter\ifx\theemail\relax    
\relax\else\vglue 5pt plus 0.02fil minus 2pt\def\\{\stdspace{\rm 
and}\stdspace} 
\cl{Email:\stdspace\tt\theemail}\fi
\ifx\theurl\relax                  
\relax\else\vglue 5pt plus 0.02fil minus 2pt\def\\{\stdspace{\rm 
and}\stdspace}
\cl{URL:\stdspace\tt\theurl}\fi\par}

\vglue 7pt plus 0.3fil minus 3pt

{\bf Abstract}
\vglue 5pt plus 0.1fil minus 2pt

\theabstract

\vglue 7pt plus 0.3fil minus 3pt

{\bf AMS Classification numbers}\quad Primary:\quad \theprimaryclass

Secondary:\quad \thesecondaryclass

\vglue 5pt plus 0.3fil minus 2pt

{\bf Keywords:}\quad \thekeywords

\vglue 10pt plus 0.5fil minus 5pt

{\small  Proposed: \theproposer\hfill Received: \receiveddate\nl
Seconded: \theseconders\hfill 
\ifx\reviseddate\relax                         
Accepted: \accepteddate                        
\else
Revised: \reviseddate                          
\fi}
\eject
}       

\let\maketitlepage\maketitlep
\let\maketitle\maketitlepage

\font\phead=cmsl9 scaled 950
\font=cmsl9 scaled 1050
\font\pnum=cmbx10 scaled 913
\font\lnum=cmbx10 
\font\pfoot=cmsl9 scaled 950
\font=cmsl9 scaled 1050
\ifplaintex
\headline{\vbox to 0pt{\vskip -4.5mm\line{\small\phead\ifnum
\count0=\startpage ISSN 1364-0380 (on line)
1465-3060 (printed) \hfill {\pnum\folio}\else\ifodd\count0\def\\{ }%
\ifx\theshorttitle\relax\thetitle\else\theshorttitle\fi\hfill{\pnum\folio}
\else\def\\{ and }{\pnum\folio}\hfill\ifx\theshortauthors\relax\theauthors
\else\theshortauthors\fi\fi\fi}\vss}}
\footline{\vbox to 0pt{\vglue 0mm\line{\small\pfoot\ifnum\count0=\startpage
\copyright\ \gtp\hfill\else
\gt, Volume \thevolumenumber\ (\thevolumeyear)\hfill\fi}\vss
}}
\else
\makeatletter
\def\@oddhead{{\small\lhead\ifnum\count0=\startpage ISSN 1364-0380 (on line)
1465-3060 (printed) \hfill {\lnum\number\count0}\else\ifodd\count0
\def\\{ }\ifx\theshorttitle\relax \thetitle \else\theshorttitle\fi\hfill
{\lnum\number\count0}\else\def\\{ and }{\lnum\number\count0}
\hfill\ifx\theshortauthors\relax 
\theauthors\else\theshortauthors\fi\fi\fi}}\def\@evenhead{\@oddhead}
\def\@oddfoot{\small\lfoot\ifnum\count0=\startpage\copyright\ \gtp\hfill\else
\gt, Volume \thevolumenumber\ (\thevolumeyear)\hfill\fi}
\def\@evenfoot{\@oddfoot}
\makeatother
\fi

\newwrite\gtoutfile
\long\gdef\makeheadfile{  
{\def\\{, }\def\s{ }
\immediate\openout\gtoutfile head.xxx
\immediate\write\gtoutfile{Proxy-for: \ifx\theasciiauthors\relax
\theauthors\else\theasciiauthors\fi\s<\ifx\theasciiemail\relax\theemail\else\theasciiemail\fi>}
\immediate\write\gtoutfile{\noexpand\\}
\immediate\write\gtoutfile{Authors: \ifx\theasciiauthors\relax
\theauthors\else\theasciiauthors\fi}
{\def\\{ }\immediate\write\gtoutfile{Title: \ifx\theasciititle\relax
\thetitle\else\theasciititle\fi}}
\immediate\write\gtoutfile{Subj-class: GT or SG or MG etc}
\immediate\write\gtoutfile{MSC-class: \theprimaryclass\ifx\thesecondaryclass\relax\else, \thesecondaryclass\fi}
\immediate\write\gtoutfile{Journal-ref: Geom. Topol. \thevolumenumber
(\thevolumeyear) \startpage-\finishpage}
\immediate\write\gtoutfile{Comments: Published by Geometry and Topology at}
\immediate\write\gtoutfile{\s\s http://www.maths.warwick.ac.uk/gt/GTVol\thevolumenumber/paper\thepapernumber.abs.html}
\immediate\write\gtoutfile{\noexpand\\}
\immediate\write\gtoutfile{}
\ifx\theasciiabstract\relax
\immediate\write\gtoutfile{\theabstract}\else
\immediate\write\gtoutfile{\theasciiabstract}\fi
\immediate\write\gtoutfile{}
\immediate\write\gtoutfile{\noexpand\\}
\immediate\write\gtoutfile{}
\immediate\closeout\gtoutfile}}  

\def\maketitlepage{\maketitlep\makeheadfile}
\let\maketitle\maketitlepage

\lognumber{471}
\received{21 June 2004}
\volumenumber{8}\papernumber{40}\volumeyear{2004}
\pagenumbers{1471}{1499}   
\revised{1 December 2004}
\published{5 December 2004}
\accepted{24 November 2004}
\proposed{Joan Birman}
\seconded{Thomas Goodwillie, Ralph Cohen}

\usepackage{amssymb, amsmath, amscd, graphicx}

\def\lbl#1{\label{#1}}

\newtheorem{theorem}{Theorem}[section]
\newtheorem{corollary}[theorem]{Corollary}
\newtheorem{lemma}[theorem]{Lemma}

\newtheorem{proposition}[theorem]{Proposition}

\newtheorem{conjecture}[theorem]{Conjecture}
\newtheorem{thm}[theorem]{Theorem}

\newtheorem{prop}[theorem]{Proposition}

\newcommand{\Out}{\operatorname{Out}}
\newcommand{\Aut}{\operatorname{Aut}}
\newcommand{\ag}{$AB$--graph }
\newcommand{\ags}{$AB$--graphs }
\newcommand{\IHX}{\text{IHX}}
\newcommand{\produ}{\pi}

\theoremstyle{definition}
\newtheorem{definition}[theorem]{Definition}

\begin{document}

\newcommand{\Ass}{\mathcal A}
\newcommand{\Com}{ \mathcal C}
\newcommand{\F}{\mathcal F}
\newcommand{\Lie}{\mathcal L}
\newcommand{\bdry}{\partial}
\newcommand{\g}{\mathcal G}
\newcommand{\fox}{d}
\newcommand{\assfox}{d^{\mathbf a}}
\newcommand{\liefox}{\mathfrak D}
\newcommand{\rpartial}{D}
\newcommand{\ls}{\mathcal{LL}} 
\newcommand{\la}{\mathcal{LA}}
\newcommand{\lc}{\mathcal{LC}} 
\newcommand{\tr}{\Tr^{\mathcal G}}
\newcommand{\Tr}{\operatorname{Tr}}
\renewcommand{\b}{\mathbf B}
\newcommand{\B}{\mathbb B}
\newcommand{\Q}{\mathbb Q}
\newcommand{\m}{\mathcal G_{AB}}
\renewcommand{\o}{o}
\newcommand{\Ge}{\mathbf G_e}
\newcommand{\GC}{\mathbf G\langle C\rangle}

\title[Morita classes in automorphism groups of free groups]{Morita classes in the homology\\of automorphism groups of free groups}
\authors{James Conant\\Karen Vogtmann}
\address{Department of Mathematics, University of 
Tennessee\\Knoxville, TN, 37996, USA}
\secondaddress{Department of Mathematics, 
Cornell Univeristy\\Ithaca, NY 14853-4201, USA}
\gtemail{\mailto{jconant@math.utk.edu}{\qua\rm 
and\qua}\mailto{vogtmann@math.cornell.edu}}

\asciiaddress{Department of Mathematics, University of 
Tennessee\\Knoxville, TN, 37996, USA\\and\\Department of Mathematics, 
Cornell Univeristy\\Ithaca, NY 14853-4201, USA}
\asciiemail{jconant@math.utk.edu, vogtmann@math.cornell.edu}

\begin{abstract} 
Using Kontsevich's identification of the homology of the Lie algebra $\ell_\infty$ with the cohomology of $\Out(F_r)$,
Morita defined a sequence of $4k$--dimensional classes $\mu_k$ in the unstable rational homology of $\Out(F_{2k+2})$.  He showed by a computer calculation that the first of these is non-trivial, so coincides with the unique non-trivial rational homology class for $\Out(F_4)$.  
Using the ``forested graph complex" introduced in \cite{exposition}, we reinterpret and generalize Morita's cycles, obtaining an unstable cycle for every connected odd-valent graph. (Morita has independently found similar generalizations of these cycles.)
The description of Morita's original cycles becomes quite simple in this interpretation, and we are able
to show that the second Morita cycle also gives a nontrivial homology
class. Finally, we view things from the point of view of a different chain complex, one which is associated
to Bestvina and Feighn's bordification of outer space. We construct cycles which appear to be the same as the Morita cycles constructed in the first part of the paper. In this setting, a further generalization becomes apparent, giving cycles for objects more general than odd-valent graphs.
Some of these cycles lie in the stable range.
We also observe that these cycles lift to cycles for $\Aut(F_r)$.
\end{abstract}

\asciiabstract{%
Using Kontsevich's identification of the homology of the Lie algebra
l_infty with the cohomology of Out(F_r), Morita defined a sequence of
4k-dimensional classes mu_k in the unstable rational homology of
Out(F_{2k+2}).  He showed by a computer calculation that the first of
these is non-trivial, so coincides with the unique non-trivial
rational homology class for Out(F_4).  Using the "forested graph
complex" introduced in [Algebr. Geom. Topol. 3 (2003) 1167--1224], we
reinterpret and generalize Morita's cycles, obtaining an unstable
cycle for every connected odd-valent graph. (Morita has independently
found similar generalizations of these cycles.)  The description of
Morita's original cycles becomes quite simple in this interpretation,
and we are able to show that the second Morita cycle also gives a
nontrivial homology class. Finally, we view things from the point of
view of a different chain complex, one which is associated to Bestvina
and Feighn's bordification of outer space. We construct cycles which
appear to be the same as the Morita cycles constructed in the first
part of the paper. In this setting, a further generalization becomes
apparent, giving cycles for objects more general than odd-valent
graphs.  Some of these cycles lie in the stable range.  We also
observe that these cycles lift to cycles for Aut(F_r).}

{\small\maketitlepage}

\section{Introduction}
The group $\Out(F_r)$ of outer automorphisms of   a free group is closely related to various classical  families of groups, including surface mapping class groups, lattices in semi-simple Lie groups and non-positively curved groups.  However, it does not actually belong to any of these families, so  theorems developed for these families do not apply directly to $\Out(F_r)$.  Although a great deal of progress has been made in recent years by adapting classical techniques to suit $\Out(F_r)$, still comparatively little is known about this group.  In particular,  invariants such as  homology and cohomology are not well understood.  The rational  homology has been completely computed only up to $n=5$ \cite{HaVo,Gerlits}, and the only  non-trivial rational homology group known prior to this paper is $H_4(\Out(F_4);\Q)\cong {\Q}$.   

A novel approach to studying these invariants was discovered by Kontsevich. He 
proved a remarkable theorem stating that $\oplus_rH^*(\Out(F_r);\Q)$  is basically the same as the homology of a certain infinite dimensional Lie algebra, $\ell_\infty$  \cite{Ko1, Ko2}. S. Morita recognized this Lie algebra as  the kernel of the bracketing map on the free Lie algebra, and used the  theory of Fox derivatives to define a family of cocycles on $\wedge\ell_\infty$ \cite{Morita}. He then applied Kontsevich's theorem to obtain a family of homology classes in $\{H_*(\Out(F_r);{\Q})\}$.  Morita was able to show that the first of these classes, which resides in $H_4(\Out(F_4);{\Q})$, is non-trivial, so in fact generates this homology group.   He  conjectured that his classes are all nontrivial,
and are basic building blocks for the rational homology of $\Out(F_r)$. In the paper \cite{Morita}, Morita was mainly interested in studying the mapping class group and he did not go into detail about  generalizing his series, but he was aware that it did generalize. We will explore such a generalization in this paper.

We use the interpretation of Kontsevich's theorem given in \cite{exposition} to translate Morita's cycles  into terms of {\it forested graphs}. This description is quite simple, and we can  easily verify that the first Morita cycle is non-trivial.  With the aid of a computer mathematics package, we are also able to show that the second Morita cycle, in $H_8(\Out(F_6);\Q)$ is non-trivial, thus supporting Morita's conjecture and producing the second known non-trivial rational homology class for $\Out(F_r)$. The graphical description of Morita's cycles also leads to a natural generalization,  giving a new cycle for every odd-valent graph.

In the last section of the paper, we diversify our portfolio of chain complexes. We introduce a new chain complex to compute $H_*(\Out(F_r);\Q)$, which we construct using the bordification of outer space introduced by Bestvina and Feighn \cite{BeFe}. 
In this new complex, we construct cycles for every odd-valent graph and conjecture that they coincide with Morita's cycles. An advantage of the new complex is that
 a further generalization becomes apparent, and we obtain a cycle associated to a more general type of graph we call an \ag\!\!. 
Although the most ambitious conjecture, that all of these cycles correspond to independent homology classes, cannot be true, it may be true for the classes coming from odd-valent graphs.  Moreover, it is possible that the classes coming from \ags\!\! generate the entire homology. 
  
Finally, we observe that all of the cycles we constructed (in the bordification context) lift to cycles in $H_*(\Aut(F_r);\Q)$, but, with the exception of a class in $H_4(\Aut(F_4);\Q)$, we do not know if they are homologically essential. A recent computer calculation of F. Gerlits \cite{Gerlits} shows that
the most ambitious conjectures here are also not true.
Namely, he found a class in $H_7(\Aut(F_5);\Q)$  which, for degree reasons, cannot arise by our construction.

\rk{Acknowledgment} Sections 2 and 3 of this paper are partly based on a set of notes by Swapneel Mahajan.
We also thank Shigeyuki Morita for an informative email in which he mentioned that he had independently found several generalizations of his original family of cycles.
The first author was partially supported by NSF grant
DMS 0305012.
The second author was partially supported by NSF grant DMS 0204185.

\section{Background and definitions}

\subsection{Kontsevich's theorem}
 In this section we briefly describe Kontsevich's Lie algebra $\ell_\infty$ and his theorem identifying the Lie algebra homology of $\ell_\infty$ with  the cohomology of the groups $\Out(F_r)$.  For a detailed exposition and proofs, we refer to  \cite{exposition}.
 
 Let $V_n$ be the $2n$--dimensional symplectic vector space with 
standard symplectic basis $p_1,\ldots, p_n,q_1,\ldots,q_n$, let $\Lie_{2n}$ be the free Lie algebra generated by
$V_n$, and let $\omega=\sum_{i=1}^n  [p_i,q_i]\in \Lie_{2n}$.  
  Kontsevich defined $\ell_n$  to be the Lie algebra consisting of those derivations of $\Lie_{2n}$ which kill $\omega$. (Recall that $f\co \Lie_{2n}\to\Lie_{2n}$ is a {\it derivation} if $f([x,y])=[f(x),y]+[x,f(y)]$).     
The bracket $[f,f^\prime]$ of two derivations is defined to be the difference 
$f\circ f^\prime-f^\prime\circ f$.
 Notice that $\ell_n\subset \ell_{n+1}$ in a natural way, and define $\ell_\infty = \underset{n\to\infty}{\lim} \ell_n$.
 
 The Lie algebra $\frak{sp}(2n)$ sits naturally as  a subalgebra of $\ell_n$, 
 and one case of Kontsevich's theorem
relates the Lie algebra homology of
$\ell_\infty$  to the homology of
the limit $\frak{sp}(2\infty)$ and to the cohomologies of the groups $\Out(F_r)$.  Specifically, we have

\begin{theorem}[Kontsevich \cite{Ko1,Ko2}]
$$PH_k(\ell_\infty) \cong  H_k(\frak{sp}(2\infty))\oplus \bigoplus_{r\geq 2} H^{2r-2-k}(\Out(F_r))$$
\end{theorem}

Here the prefix $P$ means to take primitive elements in the Hopf algebra $H_*(\ell_\infty)$.

The isomorphism in Kontsevich's theorem is defined by relating both sides to an intermediate object, called \emph{Lie graph homology}.  Briefly, the Lie graph chain complex, $\ell\g$, is spanned by finite oriented graphs, whose vertices are colored by elements of the Lie operad.  The boundary operator contracts edges one at a time, while composing the operad elements that color the endpoints of the contracting edge. As was shown in \cite{exposition}, this chain complex is equivalent to the {\it forested graph complex}, which we describe in section 4. 
The homology of $\ell_\infty$ is related to Lie graph homology via the invariant theory of $\frak{sp}(2n)$,
and the cohomology of $\Out(F_r)$ is related via the action of $\Out(F_r)$ on Outer space.

 \subsection{Fox derivatives}

Let $F_{2n}$ be the free group on $PQ_n=\{p_1,q_1,\ldots,p_n,q_n\}$.
For each $x\in PQ_n$, the classical Fox derivative is a map$${\fox_x}\co \mathbb Z[F_{2n}]\to \mathbb Z [F_{2n}],$$ where  $\mathbb Z[F_{2n}]$ is
the group ring. It is defined on $y\in PQ_n$ by the rule ${\fox_x}(y)=\delta_{xy}$, and extended to all of $\mathbb Z [F_{2n}]$ as a derivation, where the left action is multiplication and the right action
is via the augmentation $\epsilon\co \mathbb Z [F_{2n}]\to \mathbb Z $, ie,
$${\fox_x}(ab)=\left({\fox_x}a\right)\epsilon(b) + a\left( {\fox_x}b\right).$$
One can verify that this determines a well defined map.

With this starting point, we wish now to define a ``Fox derivative," on the free associative algebra,
$${\assfox_x}\co \mathcal A_{2n}\to\mathcal A_{2n},$$ where $\mathcal A_{2n}$ is the free associative
algebra on the generating set $PQ_n$. Let $x$ be an element of $PQ_n$. If $w$ is a word in the generators that does not end with
an $x$, then we define ${\assfox_x}(w)=0$; if $w=vx$, then 
${\assfox_x}(vx)=v$. 
Let $M\co \mathbb Z[F_{2n}]\to \mathcal{\widehat{A}}_{2n}$ be the Magnus embedding, where $\mathcal{\widehat{A}}_{2n}$ is the
completion to formal noncommuting power series.
 Then one can verify that the following diagram commutes:
$$
\begin{CD}
\mathbb Z[F_{2n}] @>{\fox_x}>> \mathbb Z[F_{2n}]\\
@VV{M}V @VV{M}V \\
\mathcal{ \widehat{A}}_{2n} @>{\assfox_x}>> \mathcal{\widehat{A}}_{2n}
\end{CD}
$$
It is therefore reasonable to call $\assfox_x$ an ``associative Fox derivative." 
  
The free Lie algebra embeds in the free associative algebra by sending $[a,b]$ to $ab-ba$, and 
the free associative algebra maps to the free commutative algebra $\mathcal C_{2n}$ by letting the variables commute. In this way,
we get a ``Fox derivative" ${\liefox_x}$ from the free Lie algebra to the free commutative algebra:
$${\liefox_x}\co \Lie_{2n}\to \mathcal A_{2n} \overset{{\assfox_x}}{\to}
\mathcal A_{2n} \to \mathcal C_{2n}.$$

\subsection{Morita's trace map}

Morita recognized  Kontsevich's $\ell_n$ to be  the same as a Lie algebra $\frak{h}_n$ which he had been studying in his work on mapping class groups. Let $\Lie_{2n}(k)$ be the degree $k$ part of $\Lie_{2n}$, ie,  the subspace spanned by  iterated brackets of $k$ vectors.  Now $\frak{h}_n(k)$ is defined to be the kernel of the bracketing map $V_n\otimes \Lie_{2n} (k)\to
\Lie_{2n}(k+1)$ so that we have an exact sequence:
$$0\to\ \frak{h}_n(k)\to V_n\otimes\Lie_{2n}(k)\to \Lie_{2n}(k+1)\to 0$$
We define Lie algebras $\frak{h}_n =\oplus_{k\geq 1} \frak{h}_n(k)$, and $\frak{h}_n^+= \oplus_{k\geq 2} \frak{h}_n(k)$.
For more information about $\frak{h}_n^+$, see \cite{Morita}. (In that paper, the notation 
$\frak{h}^{\mathbb Q}_{n,1}$ is the same as our $\frak{h}_n^+$.)

The isomorphism $\ell_n\to \frak{h}_n$ is given in the following way.  Given a derivation $f$ which kills $\omega$, consider the restriction $f_1$ of $f$ to $V_n$.   Then $f_1\in Hom(V,\Lie_{2n})\cong V^*\otimes \Lie_{2n}$. We use the symplectic form to identify $V^*\otimes \Lie_{2n}$ with $V\otimes \Lie_{2n}$. One easily checks that the condition that the image of $f_1$ is in the kernel of the bracket map corresponds exactly to the condition that $f(\omega)=0$.

Morita defined his trace as a map from  $\frak{h}_n$ to the free polynomial algebra $\mathcal C_{2n}$ on $2n$ variables.  Using the isomorphism above,  
  the trace becomes the following   map from $\ell_n$ to $\mathcal C_{2n}$.  
Given a derivation $f\co \Lie_{2n}\to\Lie_{2n}$, we form the Fox Jacobian
  $$\left({\liefox_{x_j}}f(x_i)\right).$$
   where $x_i,x_j\in PQ_n$
The trace of this matrix is a polynomial in $2n$ commuting variables, giving a map
$$\tau\co \ell_n\to \mathcal C_{2n},$$ 
which we will also call the   \emph{trace map}.

It will turn out that $\tau$ is a morphism of Lie algebras, provided $\mathcal C_{2n}$ is considered as an abelian Lie algebra. Indeed Morita conjectures that the image of $\tau$, together with the low degree correction, $\Lambda^3 V_n$, deriving from the first Johnson homomorphism, is exactly the abelianization of $\frak{h}_n^+$ \cite[Conjecture 6.10]{Morita}.

 \section{Graphical translation of Morita's trace}

 In this section we identify $\ell_n$ with a Lie algebra of trees,  reinterpret the trace $\tau$ as a map from this Lie algebra, and prove some of the basic properties of the trace.

We first recall how to think of generators of the free commutative, associative and Lie algebras  ($\mathcal C_{2n}, \mathcal A_{2n}$ and 
$ \Lie_{2n}$) in terms of  finite labeled trees.
 We  interpret a monomial in $\mathcal C_{2n}$  geometrically as a rooted tree with one interior vertex, whose leaves
are labeled by elements of $V_n$.  The free associative algebra
$\mathcal A_{2n}$ is spanned by rooted, labeled {\it planar} trees with one interior vertex.  The free Lie algebra $\Lie_{2n}$ is spanned by
rooted, labeled, planar  {\it  trivalent} trees.  In the Lie case (only) the trees are not linearly independent; a rooted, planar
trivalent tree with $k$ labeled leaves corresponds in a natural way to a bracket expression of $k$ letters, ie, to a generator of the
free Lie algebra.  The Jacobi relation in the free Lie algebra translates to the  IHX relation among planar trees, and the
anti-symmetry of the bracket gives the  AS   relation.

From each of the algebras $\mathcal C_{2n},\mathcal  A_{2n}$ and $\mathcal L_{2n}$, or more precisely from their underlying cyclic operads,
 we can form a Lie algebra, as in \cite{exposition}.  A generator of this Lie algebra  is again a tree of the specified type, but with no root, so that all leaves are labeled by elements of $V_n$.  We assume that each tree has at least two leaves.   We can {\it graft} a tree $X_1$ to a tree $X_2$   by identifying a leaf of $X_1$ with a leaf of $X_2$,   erasing the associated labels $v_1$ and $v_2$, and multiplying the resulting tree by a coefficient given by the symplectic product $\langle v_1,v_2\rangle$.   In the commutative and associative cases, we also contract the interior edge that we just created.
 The Lie bracket $[X_1,X_2]$ is defined by grafting $X_1$ to $X_2$ at all possible leaves, then adding up the results.   
  The Lie algebras obtained in this way from $\mathcal C_{2n}, \mathcal A_{2n}$ and $\Lie_{2n}$ are denoted $\lc_n$ $\la_n$ and $\ls_n$.   
The Lie algebra $\lc_n$ is isomorphic as a vector space to the degree $\geq 2$ part of $\mathcal C_{2n}$. 

We now wish to identify $\ls_n$ with Kontsevich's Lie algebra  $\ell_n$. Given a tree $X$ representing a generator of $\ls_n$, we  define a derivation $f_X\co\Lie_{2n}\to \Lie_{2n}$ as follows.  For each rooted tree $T$ representing a generator of $\Lie_{2n}$,  we graft $X$ onto $T$ at all possible leaves, with coefficient determined by the symplectic form, and add up the results to obtain $f_X(T)$.   The derivation $f_X$ is easily seen to kill $\omega$, so lies in $\ell_n$.
 
\begin{proposition}
Over a field of characteristic zero,  the map $\ls_n\to\ell_n$ sending $X$ to $f_X$  is an isomorphism of Lie algebras.
\end{proposition}
\begin{proof}
This follows from the fact that $\ls_n\cong \frak{h}_n$, \cite[Corollary 3.2]{levine} although the main point of that paper is that the result is  not true for the corresponding objects defined over $\mathbb Z$. The map
$\ls_n\to \frak{h}_n$ is defined on a labeled tree $X$ as a sum $\sum v\otimes X_v$, where the sum is over all univalent vertices of $X$, $v$ is the vector in $V_n$ labeling that vertex, and $X_v$ is the rooted tree formed by thinking of the chosen vertex as the root. Composing this map with the inverse of our map $\ell_n\to\frak{h}_n$, we get the indicated result.
\end{proof}

In order to understand Morita's trace map as a map from trees,
 we first need to define another type of ``partial derivative,"
$$\rpartial_{p_i}\co \ls_n\to  \mathcal L_{2n}.$$  This is given by in turn replacing each occurrence of $p_i$
by a root, and adding the results.  If $X$ is a single (unrooted) labeled trivalent tree, then 
$\rpartial_{p_i}(X)$ will be a sum of rooted labeled trivalent trees, one for each occurrence of $p_i$ at a leaf of  $X$.
 The partial $\rpartial_{q_i}$ is defined analogously. 
  
 \begin{proposition}
The trace $\tau\co\ls_n\to \mathcal C_{2n}$  is given by 
$$\tau(X)=\sum_{i=1}^n \left({\liefox_{p_i}}{\rpartial_{q_i}}(X) - {\liefox_{q_i}}{\rpartial_{p_i}}(X)\right).$$
 \end{proposition}
 
 \begin{proof}
Let $X$ be a generator of $\ls_n$, and let $f_X\co \Lie_{2n}\to\Lie_{2n}$ be the corresponding derivation, in $\ell_n$.  Then $\tau(f_X)$ is the trace of the Fox Jacobian:
$$\sum_i \left(\liefox_{p_i} f_X(p_i)+\liefox_{q_i}f_X(q_i)\right)$$
The vector $p_i$ is represented graphically as a rooted tree with a
single leaf, labeled $p_i$.  Thus $f_X(p_i)$ is a sum of rooted trees,
obtained by grafting this tree to each occurrence of $q_i$ in $X$. The
result of grafting is that the leaf formerly labeled $q_i$ becomes the
root, and we have a rooted tree representing an element of
$\Lie_{2n}$.  Thus $f_X(p_i)$ coincides with the $q_i$-th partial
derivative of $X$.  Similarly, $f_X(q_i)$ coincides with the $p_i$-th
partial derivative, but the grafting operation incurs a minus sign.
\end{proof}

We next show that the trace vanishes on all but certain very special types of trees in $ \ls_n$. 
 For any tree $T$, we define the {\it interior} of $T$  to be the tree obtained by  removing all leaves and the edges incident to those leaves. 
 If the
interior of
$T$ has no trivalent vertices (ie is a line), we say $T$ is {\it linear}.  If the interior of $T$ has an odd number of edges, we say $T$ is {\it odd};
otherwise
$T$ is {\it even}.

\begin{lemma}
Let $T$ be a generator of $\Lie_{2n}$, represented by a rooted, labeled, planar trivalent tree.  Then $\liefox_{p_i}(T)=0 $ unless $T$ is a
linear tree with exactly one leaf labeled $p_i$ at maximal distance from the root.  If the other leaves
 are labeled $a_1,\ldots,a_k$, then 
$\liefox_{p_i}(T)=\pm a_1\ldots a_k.$ 
The sign is determined by the fact that
 $\liefox_{p_i}[a_1,[a_2,\ldots,[a_k,p_i]\ldots]]=a_1\ldots a_k$.
\end{lemma}

\begin{proof}
If $T$   is of the specified form, then consider the unique geodesic path from the root to the leaf labeled $p_i$.  We can use AS relations, if necessary, to put all of the leaves $a_j$ onto one side of this  geodesic, so that  we may assume $T$ is the tree corresponding to   $$[a_1,[a_2,\ldots[a_k,p_i]\ldots]].$$  Consider the image $P$ of $[a_1,[a_2,\ldots[a_k,p_i]\ldots]]$ in the free associative algebra $\mathcal A_{2n}$.  If no $a_j$ is equal to $p_i$, then the only term of $P$ which ends in $p_i$ is $a_1a_2\ldots a_kp_i$, so that $\liefox_{p_i}(W)=  a_1\ldots a_k.$  On the other hand, suppose some $a_j$ is equal to $p_i$, with $j<k$ minimal. Let $\beta=[a_{j+1},\ldots,[a_k,p_i]\cdots]$.
 Then the only terms that could possibly end in $p_i$ come from 
 \begin{align*}
 a_1\cdots a_{j-1}[p_i,\beta]=a_1\cdots a_{j-1} p_i\beta-  a_1\cdots a_{j-1}\beta p_i.
 \end{align*}
By induction on the number of occurrences of $p_i$, we know that $\liefox_{p_i} \beta =a_{j+1}\cdots a_k$, so that $\liefox(a_1\cdots a_{j-1} p_i\beta)=a_1\cdots a_{j-1} p_i\liefox(\beta)= a_1\cdots a_k$.  We also have
 $$\liefox_{p_i} ( a_1\cdots a_{j-1}\beta p_i)= a_1\cdots a_{j-1}\beta,$$ but this second term vanishes when we pass to $\mathcal C_{2n}$ since $\beta$ is a commutator.
 
 If $T$ has any form other than the one specified, then its image $P$ in $\mathcal A_{2n}$ can be grouped  into terms so that $\liefox_{p_i}(T)$ contains a commutator as a factor in each term, which becomes zero
when we pass to $\mathcal C_{2n}$.
\end{proof}

The tree $T$ corresponding to $[a_1,[a_2,\ldots,[a_k,p_i]\ldots]]$ is shown in Figure~\ref{dnv}a.
\begin{figure}
\begin{center}
$\underset{\text{\small(a)}}{\includegraphics[width=1.2in]{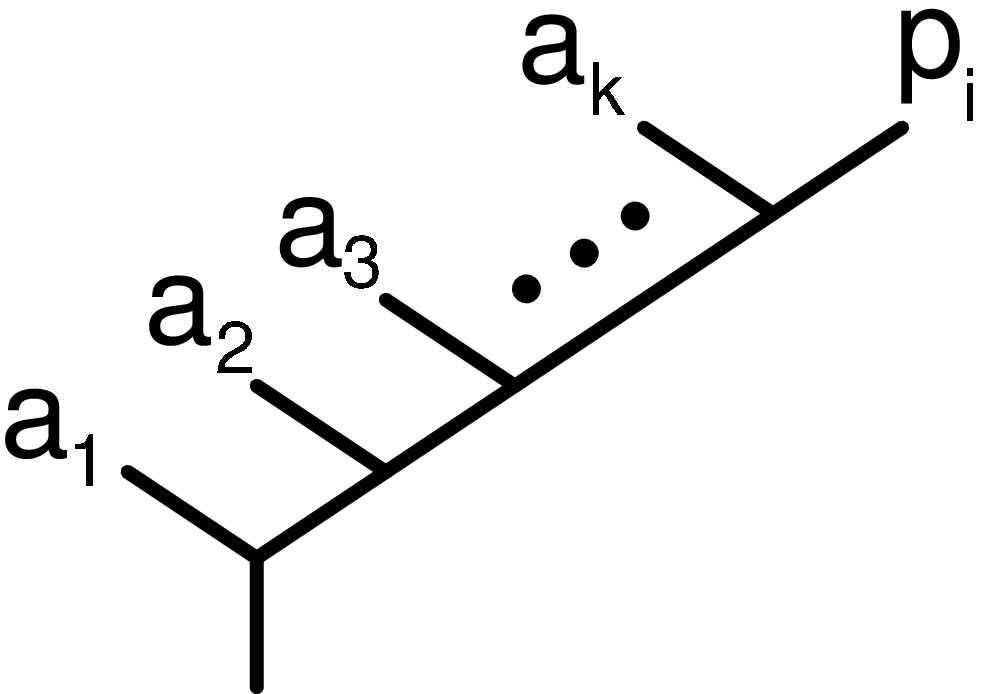}}\hspace{1in}\underset{\text{\small(b)}}{\includegraphics[width=1.5in]{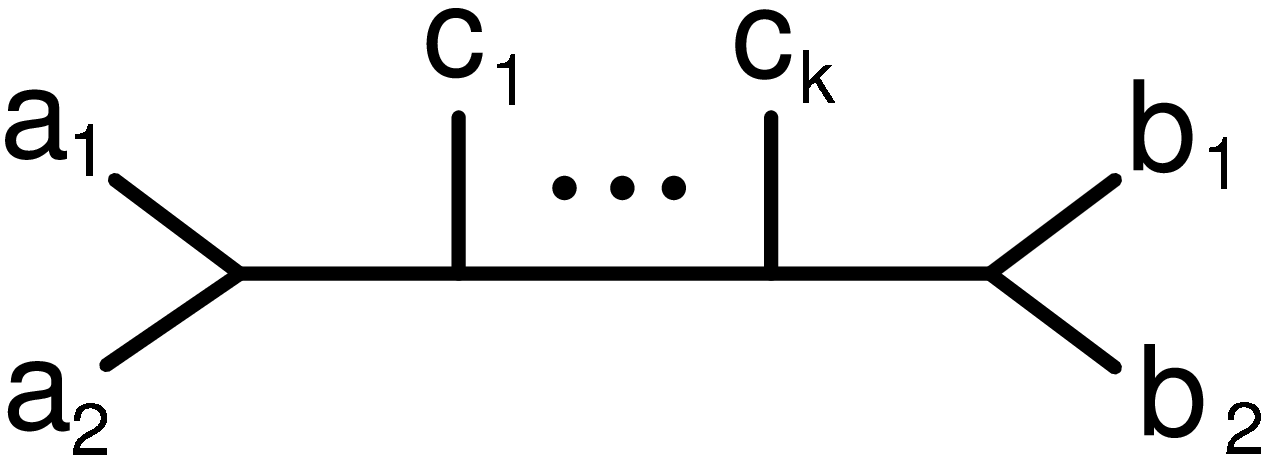}}$
\caption{(a) The unique tree on which $\liefox_{p_i}$ does not vanish, and (b)
an element of $\ls_n(k+2)$.
}\lbl{dnv}
\end{center}
\end{figure} 
 Note that this is
$(-1)^k[[\ldots[p_i,a_k],\ldots,a_2],a_1]$, which corresponds to the tree obtained by flipping all of the leaves labeled $a_i$ across
the axis from the root to $p_i$. This observation allows us to easily compute the trace of any generator of $\ls_n$.

\begin{corollary}\label{a2}
Let $X$ be a nonzero generator of $\ls_n$, represented by a labeled planar trivalent tree.  Then $\tau(X)=0$ unless $X$ is of the form
shown in Figure~\ref{dnv}b, with $k$ odd.  In this case,
$$\tau(X)=2c_1\ldots c_k\big(-\langle a_1,b_1\rangle a_2b_2 +\langle a_1,b_2\rangle a_2b_1+\langle a_2,b_1\rangle a_1b_2-\langle
a_2,b_2\rangle a_1b_1\big),$$
where $\langle\cdot\, ,\cdot \rangle$ is the symplectic form.
\end{corollary}

Note that  $a_1\neq a_2$ and $b_1\neq b_2$, since $X$ is non-zero.  If $a_1=p_i$, $b_1=q_i$ and $a_2$ and $b_2$ are not
paired, then $\tau(X)=-2a_2c_1\ldots c_kb_2$.

\begin{lemma}~\label{lem35}
For any $X,Y\in \ls_n$, 
$\tau([X,Y])=0.$
\end{lemma}

\begin{proof}
Each term in $[X,Y]$ is obtained by grafting a leaf of $X$ to a leaf of $Y$ to obtain a new tree $Z$.  If either $X$ or $Y$ is
non-linear, then so is $Z$.  If $X$ and $Y$ are both linear, and both even or both odd, then the interior of $Z$ is either non-linear,
or is odd linear, so has trace zero.  The only possible way to obtain a term with a non-zero trace is if $X$ is odd linear and $Y$ is
even linear, or vice-versa, and we are grafting extremal leaves of $X$ and $Y,$ say labeled $p_i$ and $q_i$.

In order for the trace of $Z$ to be non-zero, it must have an extremal pair of leaves labeled $q_j$ and $p_j$ for some $j$.
 In this case, there is another term $Z'$ in the bracket, obtained by grafting these two leaves.
We claim that the  trace of $Z$ cancels with the trace of $Z'$.  

In the following figure, $X$ is on the left, $Y$ is on the right.
\medskip

\centerline{\includegraphics[width=1.7in]{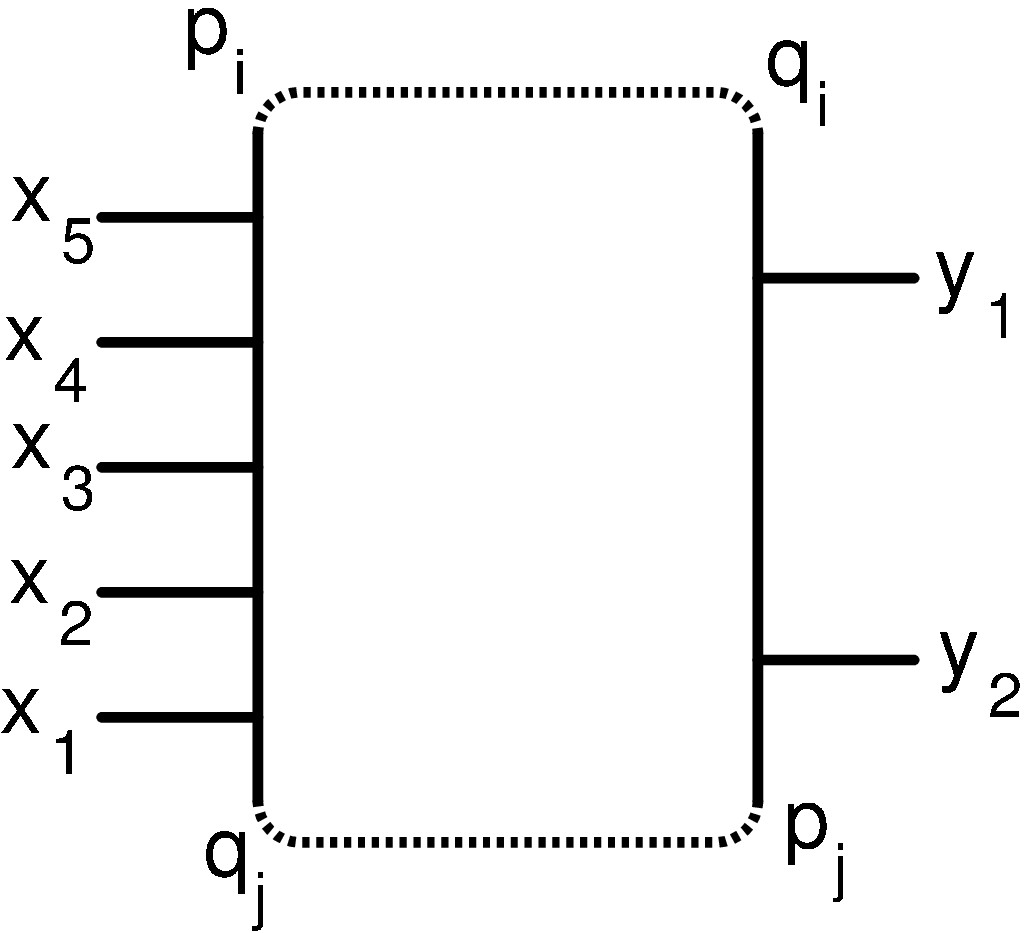}}

$Z$ is obtained by joining the lower two leaves, whereas $Z^\prime$ is obtained by joining the upper two leaves and multiplying by $-1=\langle q_j,p_j\rangle$.
Let $\produ=2x_1\cdots x_sy_1\cdots y_t$, where the $x_k$'s and the $y_k$'s are the labels of the leaves of $X$ and $Y$ which are present in both $Z$ and $Z^\prime$. 
There are an odd number of such leaves, and $\produ$ is the trace, up to sign.
We can draw $Z$ as a horizontal line from $q_i$ to $p_i$ with all of the other edges pointing upwards from this line. By Corollary~\ref{a2} the trace of this is $\produ$. We can draw $-Z^\prime$ as a horizontal edge from $q_j$ to $p_j$ with the other edges pointing upwards. The trace of $-Z^\prime$ is therefore also $\produ$ which implies that $\tau(Z^\prime+Z)=-\produ+\produ=0$ as desired. 
\end{proof}

Lemma~\ref{lem35} shows that $\tau$ is a Lie algebra homomorphism, if we give $\mathcal C_{2n}$
an \emph{abelian} Lie algebra structure.

We can extend the trace  to a function on the exterior algebras $\wedge \ls_n\to \wedge \mathcal C_{2n}$, by defining
$$\Tr(X_1\wedge\ldots\wedge X_k)=\tau(X_1)\wedge\ldots\wedge\tau(X_k).$$  Recall that the Lie algebra boundary map is given by
$$\bdry(X_1\wedge\ldots\wedge X_k)=\sum_{i<j} (-1)^{i+j} [X_i,X_j]\wedge X_1\wedge\ldots\wedge\hat X_i\wedge \ldots\wedge\hat X_j\wedge\ldots\wedge X_k$$
Then the following is an immediate consequence of Lemma~\ref{lem35}:

\begin{lemma}
Trace is a cocycle, ie, $\Tr(\partial X)=0$ for any $X\in \wedge \ls_n$.  
\end{lemma}

\section{Forested graphs}\label{fg}

In this section we define the {\it graphical trace} as a map between chain complexes of graphs.  The domain of this trace is the {\it forested graph complex}, which computes the rational cohomology of $\Out(F_r)$.  We explain its precise relationship with the map $Tr$ defined in the previous section, and show that it is a cocycle, both by exploiting this connection and by giving an independent combinatorial proof.  

Recall from \cite{Ko1,Ko2} (see also \cite{exposition}) that the Lie algebra  $\ell_n$ contains the symplectic
Lie algebra
$\frak{sp}(2n)$ as a subalgebra, and that the Lie algebra homology  of $\ell_n$ can be computed as the homology of the subcomplex of
$\frak{sp}(2n)$--invariants in the exterior algebra $\wedge\ell_n$.  Under the identification of $\ell_n$ with $\ls_n$, $\frak{sp}(2n)$ is spanned by trees with exactly two leaves.  (Such trees also span a copy of $\frak{sp}(2n)$ in $\la_n$ and in $\lc_n$.) 

The $\frak{sp}(2n)$ invariants in  $\wedge \ell_n$ can be described as forested graphs 
modulo IHX relations \cite{exposition}.   Here a {\it forested graph} is a pair $(G,F)$, where $G$ is a finite graph all of whose vertices are either trivalent or bivalent, and  $F$ is a (possibly disconnected) acyclic subgraph containing all of the vertices of $G$.  
A forested graph comes with an  {\it orientation}, given by ordering  the edges of $F$ up to even permutations. The {\it forested graph complex} $f\g$ is the chain complex spanned by all forested graphs modulo IHX relations.  The boundary of $(G,F)$ is the sum, over all $e\in G-F$ for which $F\cup e$ is still a forest, of the forested graphs $(G, F\cup e)$.
 
There is a natural chain map $\psi_n\co \wedge \ls_n\to f\g$, defined in the following way.
A generator of $\wedge  \ls_n$ consists of a wedge $X_1\wedge\ldots\wedge X_k$  of several trees whose leaves are
labeled by elements of $V_n$.  A \emph{pairing} $\pi$ of these leaves is a partition of the leaves into two-element subsets; this also gives a pairing of the associated labels.  For each possible pairing $\pi$, glue the paired leaves together to form a trivalent graph $G_\pi$, and  use the interiors of the $X_i$ to form a forest $F_\pi$ in $G_\pi$.  The edges of the interior of each $X_i$ come with a natural ordering (see \cite{exposition}), and the ordering of the $X_i$ then gives an ordering of all edges of $F_\pi$.   Now
take the symplectic product of each pair of labels, and multiply them together to form a coefficient $w(\pi)$.   Then $$\psi_n(X_1\wedge\ldots\wedge X_k)=\sum_{\pi}w(\pi)(G_\pi,F_\pi).$$  It follows from results in section 2.5 of \cite{exposition} that the limit map $\psi_\infty\co\wedge \ls_\infty\to f\g$ is an isomorphism on homology.

There is a similar chain map  $\psi_\infty$  from $\wedge\lc_\infty$ to a complex of graphs (and one for $\la_n$).  For $\lc_\infty$, the relevant graph complex is  called $\widehat{c\g}$.  A generator of $\widehat{c\g}$ is a finite {\it oriented} graph $G$, which may have vertices of any valence, including $1$.  An orientation on $G$ is determined by ordering the vertices and orienting each edge of $G$. If all vertices of $G$ have odd valence, this is equivalent to giving an ordering of the half-edges coming into each vertex. The map $\psi_n$ is again a sum, over all pairings, of the graphs $G$ obtained by gluing the leaves according to the pairing, and multiplying the result by a coefficient coming from the symplectic products of the leaf labels.  Since each $X_i$ corresponds to a vertex of $G$, the ordering of the $X_i$ gives an ordering of the vertices of $G$, and the sign of the symplectic product of the labels gives an orientation on each edge. 

Recall that the part of $\mathcal C_{2n}$ of degree $\geq 2$, which we denote $\mathcal C^+_{2n}$, is isomorphic as vector space to $\lc_n$.
Thus we have  maps $\psi_n\co \wedge \mathcal C^+_{2n}\cong \wedge\lc_n\to \widehat{c\g}$. 
Extend $\psi_n$ to all of $\wedge\Com_{2n}$ by considering a linear generator as a tree with one basepoint and one labeled vertex.  Graphs in the image of $\psi_n$  may then have univalent vertices.

We are finally ready to define the graphical version of the trace map.  
\begin{definition}
The {\it graphical trace map}
$\tr\co f\g\to \widehat{c\g}$ is defined in the following way. Let $(G,F)$ be a forested graph. Then $\tr(G,F)=0$ unless $F$ consists of a disjoint union of linear trees $T_1,\ldots, T_k$, each with an even
number of edges, and for each $T_i$ there is an edge $e_i$ of $G-F$ joining the two ends of $T_i$.  In this case,  $\tr(G,F)$ is the graph 
 obtained by collapsing each cycle $T_i\cup e_i$ to a point, multiplied by $(-2)^k$.  The
orientation on  $\tr(G,F)$ is given by the induced ordering of the edges coming into each vertex (note that all vertices of $\tr(G,F)$ have odd valence, and may even have valence $1$).
\end{definition}

\begin{theorem}\label{comdiag} The following diagram is commutative.
$$\begin{CD}
\wedge \ls_n @>{\Tr}>> \wedge \Com_{2n}\\
@VV{\psi_n}V @VV{\psi_n}V\\
f\g @>{\tr}>> \widehat{c\g}
\end{CD}
$$
\end{theorem}
\begin{proof}
Let $X\in \wedge \ls_n$. 
The only  way that $\Tr(X)$ can be nonzero is if  $X$ is a wedge of  linear trees $X_i$, each with ends labeled by some matched pair $p_j$ and  $q_j$. Applying $\psi_n$  to $\Tr(X)$ glues up the rest of the leaves in some way.
If we apply $\psi_n$ first to $X$, then any pairing which does not match the ends of each $X_i$ is sent to zero by $\tr$.  If the ends of each $X_i$ are matched by a pairing $\pi$,  then applying $\tr$ to $(G_\pi,F_\pi)$  gives a graph where the rest of the leaves
have been glued together in some way. 
\end{proof}

\begin{corollary} $\tr\co f\g \to \widehat{c\g}$ is a cocycle.
\end{corollary}
\begin{proof} This follows since $\Tr$ is a cocycle and the left-hand vertical map is surjective.\end{proof}

\begin{prop}
 Any homology class detected by $\Tr$ is also detected by $\tr$. 
\end{prop} 
\begin{proof}
The subcomplex of $\mathfrak{sp}(2n)$--invariants of $\wedge \ls_n$ is quasi-isomorphic to the whole complex. But $\Tr ((\wedge \ls_n)^{\mathfrak{sp}(2n)})\subset  (\wedge \Com_{2n})^{\mathfrak{sp}(2n)}.$ Now it suffices to observe that $\psi_n$ is an isomorphism from the spaces of $\mathfrak{sp}(2n)$--invariants
onto $f\g$ and $\widehat{c\g}$ in the limit when $n\to\infty$.
 \end{proof}

Recall that it is the primitive part of $H_*(\ell_\infty)$ that contains the cohomology of $\Out(F_r)$. 
On the level of forested graphs, this  is the homology of the subcomplex spanned by connected graphs. It is shown in \cite{exposition} that the quotient of this subcomplex by the subspace of graphs with separating edges is quasi-isomorphic to the whole subcomplex, and that the bivalent vertices contribute to the $\frak{sp}(2\infty)$ part of the homology, and not to the $\Out(F_r)$ part.  In fact, if we let $f\g^\prime$ be the quotient of $f\g$ by the subspace spanned by graphs with separating edges,  graphs with bivalent vertices and  disconnected graphs, we have
$$H_*(f\g^\prime)\cong \underset{r}{\oplus} H^*(\Out(F_r);\Q).$$
Since we are primarily interested in the cohomology of $\Out(F_r)$, we would like the graphical trace to have domain $f\g'$.  To this end, we define $c\g^\prime$ to be the quotient of the space of oriented connected graphs with vertices of odd valence modulo the subspace spanned by graphs with separating edges. 
Note that a graph with no separating edges also has no univalent vertices, so $c\g^\prime$ is generated by graphs with odd valence at least 3. We then have
the following result.

\begin{proposition} The trace map induces a cocycle  $\tr\co f\g^{\prime}\to c\g'$. 
\end{proposition}
\begin{proof}
The trace map  preserves the number of connected components, and vanishes on graphs with bivalent vertices. Furthermore, the image lands in the subspace of graphs with vertices of odd valence. Finally the trace of a graph with separating edges is also a graph with separating edges.
\end{proof}

The fact that $\tr$ is a cocycle also has a direct, combinatorial argument, independent of the trace's origin in the world of Fox derivatives:

\begin{proposition} For any forested graph $(G,F)$, $\tr(\bdry (G,F))=0$.
\end{proposition}

\begin{proof}

The boundary $\bdry(G,F)$ is the sum, over all possible ways of adding an edge $e$ to $F$, of $(G, F\cup e)$.  If any vertex of $F$ is
trivalent (so that
$F$ is not the union of a set of linear trees), then the same will be true of $F\cup e$, so $\Tr(G, F\cup e)=0$.  So we may assume all
trees in 
$F$ are linear. If $e$ does not join two endpoints of trees in $F$, then again there will be a trivalent vertex in $F\cup e$, so
$\Tr(G,F\cup e)=0$. If $e$ joins an endpoint of a tree $T$ with an endpoint of another tree $T'$, then the tree 
$T\cup e\cup T'$ has $e(T)+e(T')+1$ edges.  In order for the term $\Tr(G,F\cup e)$ of $\Tr(\bdry(G,F))$ to be nonzero, we need $e(T) +
e(T')+1$  even,  and we need there to be another edge $f\neq e$ of $G-F$ joining the opposite ends of $T$ and $T'$.  In this case,
$\Tr(G,F\cup e)$ cancels with
$\Tr(G,F\cup f)$,  which is also a term in $\bdry(G,F)$, so that again $\Tr(\bdry(G,F))=0$.\end{proof}

\section{Non-triviality of the first two Morita cycles}

\begin{figure}
\begin{center}
\includegraphics[width=1in]{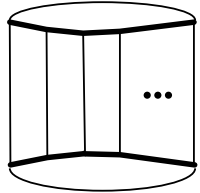}
\caption{The forested graph $(G_k, F_k)$. There are $k$ vertical edges, and the darker edges represent $F_k$.}\label{thetak}
\end{center}
\end{figure}

Morita's original trace maps \cite{Mtrace} correspond to the
compositions $\mu_k$ of our $\tr$ with projection onto the linear
subspace of $c\g^\prime$ spanned by the graph $\Theta_k$ for $k$ odd,
where $\Theta_k$ has two vertices connected by $k$ edges.  The graph
$\Theta_k$ is the trace of the forested graph $\frac{1}{4}(G_k,F_k)$
shown in Figure~\ref{thetak}, so the projection, and hence the cocycle
$\mu_k$, is non-trivial.  Appealing to the identification of
$\underset{r}{\oplus} H^*(\Out(F_r);\Q)$ with $H_*(f\g^\prime)$, it is
easy to check that the {\it cocycle\/} $\mu_k$ corresponds to a {\it
cycle\/} in $H_{4k}(\Out(F_{2k+2});\Q)$. In this section, we address
the question of whether the cocycles $\mu_k$ determine nontrivial
cohomology classes in $H^*(f\g^\prime)$, and hence determine
nontrivial homology classes in $H_*(\Out(F_r);\Q)$.

In order to prove that $\mu_k$ is nontrivial on the level of cohomology, it suffices to produce a cycle $z$ in $f\g^\prime$ with $\mu_k(z)\neq 0$.  Since we know that $\mu_k(G_k,F_k)\neq 0$,  
we calculate the boundary of $(G_k,F_k)$:
$$\bdry(G_k,F_k)=(G_k,T_1) + (G_k,T_2)+\ldots+(G_k,T_{k}),$$
where $T_i$ is the maximal tree consisting of $F_k$ together with the $i$th edge connecting the two trees of $F_k$.  It is not obvious whether this is equal to $0$ in $f\g^\prime$, ie, whether it lies in the IHX subspace of the vector space $t\g$ spanned by pairs $(G,T)$ with $G$ trivalent and $T$ a maximal tree.  
Our strategy will be to show that the subspace spanned by forested graphs with zero trace maps onto $t\g/\IHX$, so that the forested graph $(G_k,F_k)$ can be adjusted by adding traceless forested graphs to become a cycle with non-zero trace.

We begin by noting  two elementary lemmas which reduce the number of pairs $(G,T)$ we need to generate $t\g/\IHX$.

\begin{lemma}
If an edge of $G-T$ joins the endpoints of an edge of $T$, then $(G,T)=0\in f\g^\prime$.  
\end{lemma}
\begin{proof}
Apply an IHX relation to the edge of $T$ joined by the edge of $G-T$ to get a relation of the form $(G^\prime, T^\prime)=2(G,T)$, where $G^\prime$ is a graph that has a separating edge.
\end{proof}

\begin{lemma}
Any forested graph $(G,T)\in t\g$ is equivalent  to a sum of graphs $(G_i,L_i)$, where
$L_i$ has no trivalent vertices (ie it is a line), and an edge of $G_i-L_i$ joins the ends of $L_i$.
\end{lemma}
\begin{proof}
Fix an edge of $G-T$, and choose a geodesic in $T$ which connects the endpoints of the edge. If this geodesic equals $T$, then we are done. Otherwise, choose an edge of $T$ which is incident to the geodesic, but not in it, and apply an IHX relation to that edge. The result is a difference of two graphs $(G^\prime, T^\prime)$ and $(G^{\prime\prime},T^{\prime\prime})$ in which the geodesic has grown in length.  Iterate until the geodesic fills the tree.
\end{proof}

We will call a pair $(G,L)$, with $L$ a line whose endpoints are joined by an edge of $G-L$, a {\it chord diagram}, and an edge of $G-L$ which joins successive vertices of $L$ an {\it isolated chord}.   Thus $t\g/\IHX$ is generated by chord diagrams with no isolated chords.

\begin{proposition}
The cocycle $\mu_3$  represents a non-trivial class in the cohomology of the forested graph complex.
\end{proposition}

\begin{proof}  We need to find a cycle $z$ with $\mu_3(z)\neq 0$.  We know that $\mu_3(G_3,F_3)\neq 0$, so we calculate
\begin{align*}
\bdry(G_3,F_3)&=(G_3,T_1) + (G_3,T_2)+(G_3,T_3)\\
&=2(G_3,T_1) + (G_3,T_2).
\end{align*}

The graph $G_3$ has fundamental group $F_4$, and the quotient $t\g/\IHX$ is spanned by chord diagrams $(G,L)$ with no isolated chords.  
In rank 4,  there are only four  chord diagrams with no isolated chords.  If the vertices of $L$ are labeled $1,2,3,4,5,6$, then the four chord diagrams are:

\vspace{.5em}

\noindent$\underset{A= (13)(25)(46)}{\includegraphics[width=1in]{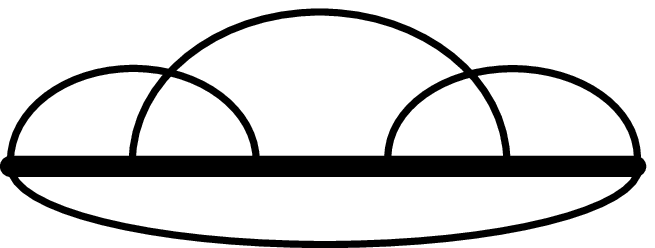}}$\hfill
 $\underset{B= (14)(25)(36)}{\includegraphics[width=1in]{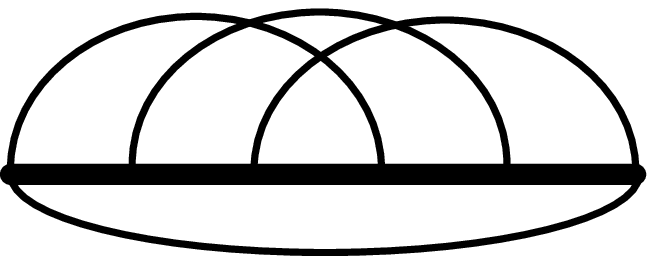}}$\hfill
 $\underset{C=(14)(26)(35)}{\includegraphics[width=1in]{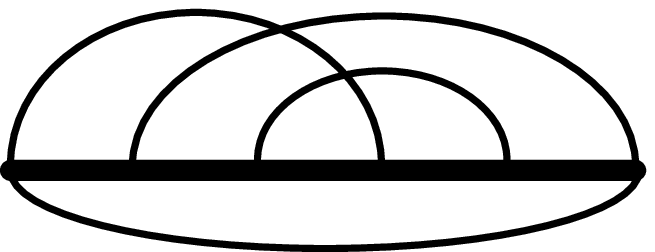}}$\hfill
 $\underset{D=(16)(24)(35)}{\includegraphics[width=1in]{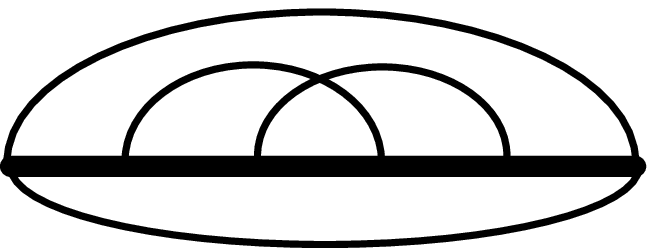}}$

\vspace{.7em}

An IHX relation on the edge (12) of  $L$ in $D$ shows $D=-2C$.  Two IHX relations, using edges $(45)$ of $L$ in $B$ and then the image of $(56)$ give $B=-2C$. 
Thus $t\g/\IHX$ has dimension at most two, and is spanned by $A$ and $C$. 

Consider the following two subforests $F^\prime$ and $F^{\prime\prime}$ of $G_3$.

\vspace{.5em}

\centerline{$\underset{(G_3,F^\prime)}{\includegraphics[width=1.5in]{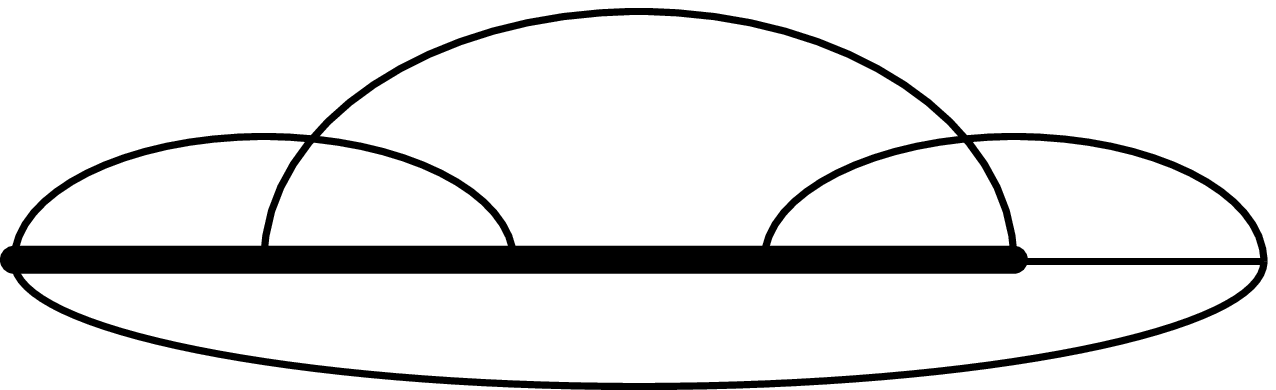}}$\hspace{3em}
 $\underset{(G_3,F^{\prime\prime})}{\includegraphics[width=1.5in]{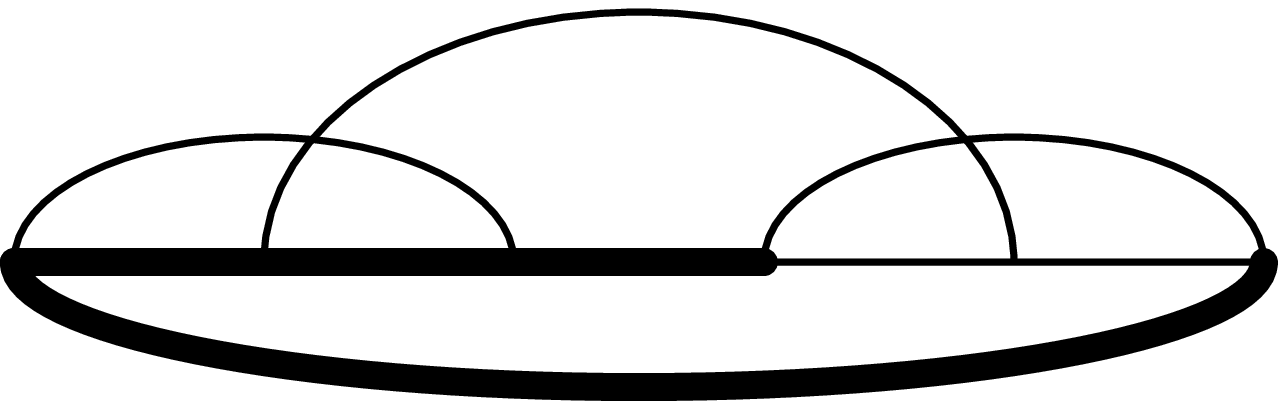}}$}
 
\vspace{.7em}

We compute
\begin{alignat*}
\bdry(G_3,F^\prime)&= C &\hbox{ mod } \IHX,&&\\
\bdry(G_3,F^{\prime\prime})&= 4C-A &\hbox{\qquad  mod } \IHX.&&
\end{alignat*}
Both $(G_3,F^\prime)$ and $(G_3,F^{\prime\prime})$ have zero trace. Thus the image of the subspace of traceless graphs under the boundary map is all of $t\g/\IHX$, so that traceless graphs can be added to $(G_3,F_3)$ to make a cycle $z$, with non-zero trace $\mu_3(z)$.

In fact, for $r=3$ it is easy to write down a cycle $z$ explicitly.  We compute that modulo IHX, the boundary of $(G_3, F_3)$ is equal to $-3A$, so that $z=(G_3,F_3)-3(G_3,F^{\prime\prime})+12(G_3,F^\prime)$ is a cycle with non-zero trace.
\end{proof}
 
For $k$ odd and bigger than $3$, the calculations become much too cumbersome to do by hand, so we  have written a computer program that does them.  We can use this to show that the second Morita class is non-trivial:

\begin{proposition}
The cocycle $\mu_5$,  corresponding to the second Morita cycle, represents a non-trivial cohomology class in $H^*(f\g^\prime)$, and hence a nontrivial homology class in $H_8(\Out(F_6);\Q)$.
\end{proposition}
\begin{proof}

Up to reflectional symmetry, there are $184$ chord diagrams  for $k=5$ with no isolated chords.
These are not linearly independent. In fact one gets a relation by choosing any chord and sliding its
endpoints to the end of the chain, modulo IHX. 
Computer calculations give  148 independent sliding relations, so 
$$\hbox{dim }t\g/IHX \leq 36.$$
As we did in the case $k=3$, we now try to kill this $36$--dimensional vector space with boundaries of traceless elements. This is accomplished by looking at the two types of forested graphs in Figure ~\ref{classill}.
\begin{figure}
$\underset{\text{Type One}}{\includegraphics[width=2in]{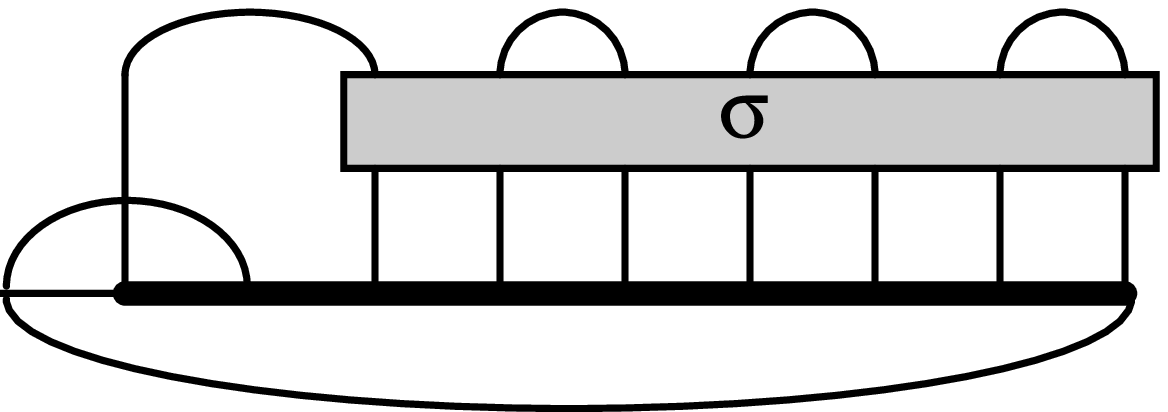}}\hfill\underset{\text{Type Two}}{\includegraphics[width=2in]{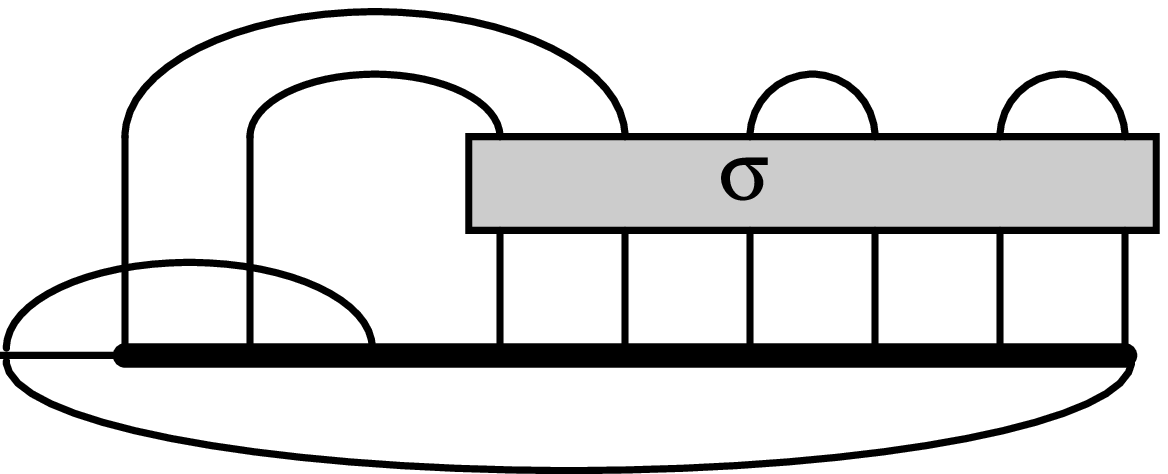}}$
\caption{Two types of graphs. $\sigma$ is a permutation.}\lbl{classill}
\end{figure}
These are formed by taking a chord diagram where the left most chord has one or two other chord feet between it. 
Then one removes the left-most edge from the chain.
Call these types of forested graphs ``type one" and ``type two," respectively. There are 39 forested graphs of type one, and $32$ forested graphs of type two.
 First note that
the trace of each of these is zero.
 The boundary of
each  has three terms, since there are only three possible edges that can be added to the chain to make a
tree.  These three terms are shown for a type one graph in Figure~\ref{typeone}.
 \begin{figure}
 $\includegraphics[width=1.5in]{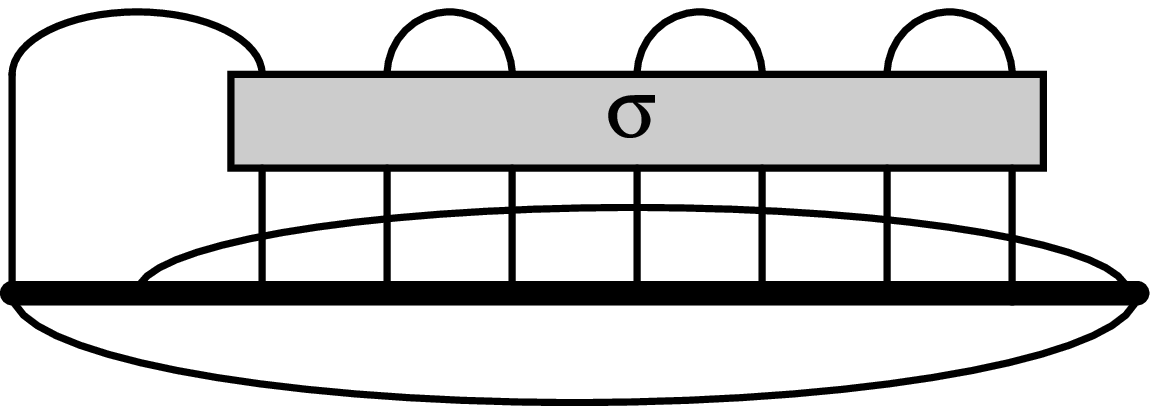}\overset{+}{\phantom{M^S}}\includegraphics[width=1.5in]{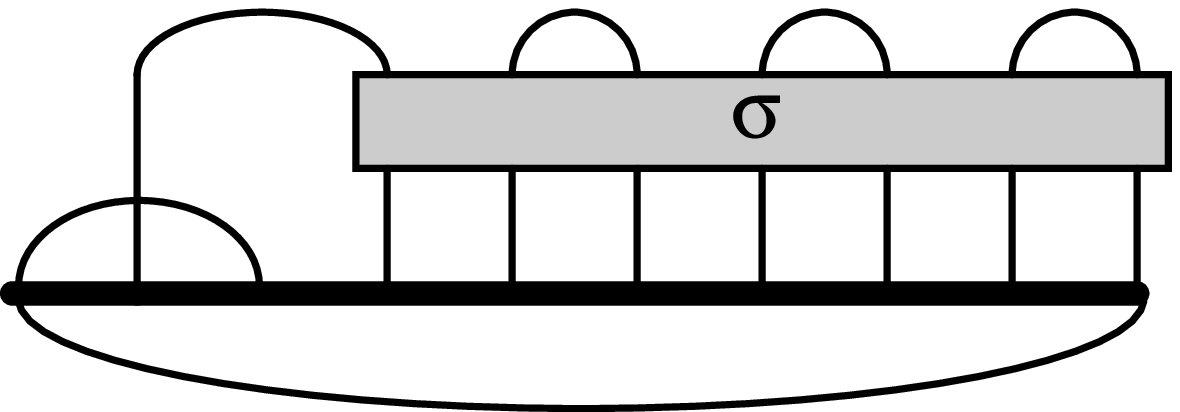}\overset{+}{\phantom{M^S}}\includegraphics[width=1.5in]{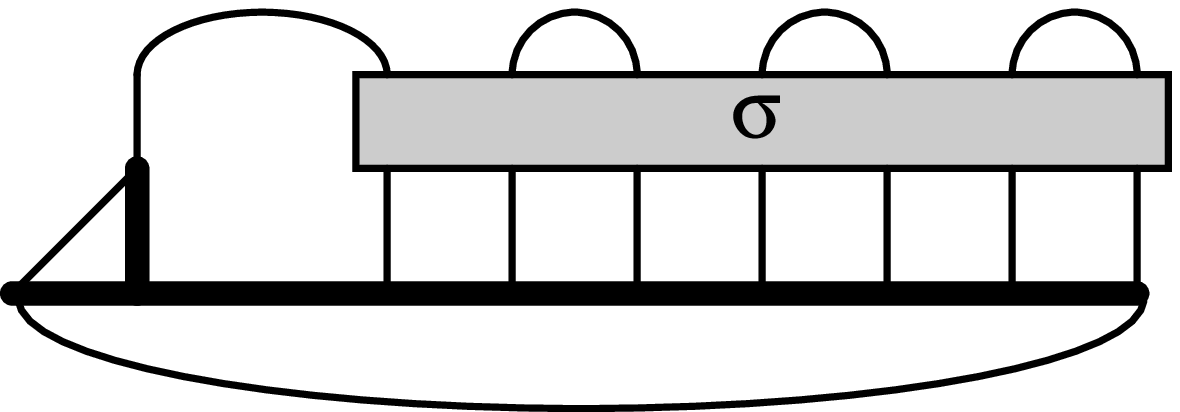}$
 \caption{The boundary of a type one graph.}\label{typeone}
 \end{figure}
 
The first two terms are already chord diagrams, and are of the form $X+\rho(X)$, where $\rho(X)$ is the diagram where the left endpoint of the first chord gets rotated around to
 be the last endpoint on the chain. 
A single IHX relation applied to the vertical edge of the tree 
will turn the last term of Figure~\ref{typeone} into a sum of two terms, one of which is equal to $X$, and the other of which is zero
because it contains an isolated chord.  
 Thus, the boundary of a type one forested graph is of the form
$2X + \rho(X)$.  All 39 of these boundaries are independent.

The boundary of a type two forested graph is of the form
$Y+\rho(Y)+Y_{(34)}-Y_{(243)}-Y_{(23)}$. Here the notation $Y_\sigma$ where $\sigma\in \Sigma_{10}$
means the graph formed from $Y$ by re-gluing the endpoints of the chords according to the permutation $\sigma$.  All 32 of these are independent, and together with the type one boundaries they give 71 independent vectors which, together with the sliding relations, span all of $t\g$. In more detail, if you take the 36-dimensional vector space and further mod out by type two relations, a $7$ dimensional vector space results. Further modding out by type one relations yields a vector space of dimension $0$.
 
 This shows that the boundary map is onto, even when restricted to
chains of trace zero. Thus we can conclude that the $k=5$ Morita class is nontrivial,
since we can take any chain on which the cocycle evaluates nontrivially, and add traceless chains to it to make a
cycle.
 \end{proof}

Recall that the virtual cohomological dimension  of $\Out(F_6)$ is equal to $9$.  Our computations show that this top-dimensional homology for $\Out(F_6)$ vanishes:

\begin{corollary}
$H_9(\Out(F_6);\Q)$ is trivial. 
\end{corollary}
\begin{proof}
The end of the forested graph complex in rank 6 is $f\g^\prime_2\to f\g^\prime_1\to 0$, where the subscript indicates the number of trees in a forest.  In the proof of Proposition 5.4, we showed the boundary map $ f\g^\prime_2\to f\g^\prime_1$ is onto, even when restricted to traceless chains.  The group $f\g^\prime_1/\operatorname{im}(\bdry)$  corresponds to the top cohomology $H^9(\Out(F_6);\Q)\cong H_9(\Out(F_6);\Q)$.
\end{proof}

\section{Generalized Morita cycles}

\subsection{Cocycles parameterized by odd-valent graphs}
To obtain the Morita traces, we composed $\tr\colon f\g'\to c\g'$ with projection onto the subspace of $c\g'$ spanned by $\Theta_k$.  In fact, projection onto the span of any (odd-valent) generator of $c\g'$  gives a cocycle, which we may think of as having values in ${\mathbb Q}$.  

 \begin{proposition}
 Let $\Theta\in c\g^\prime$ be nontrivial, and let $\mu_\Theta$ be the composition of $\tr$ with projection onto the span of  
$\Theta$.  Then $\mu_\Theta$ is nontrivial as a cocycle.
\end{proposition}

\begin{proof}  Given an odd-valent oriented graph $\Theta$, we need to construct a forested graph $(G,F)$ with $\mu_\Theta(G,F)=(-2)^v\Theta$, where $v$ is the number of vertices of $\theta$.  To do this, we let the orientation on $\Theta$ be given as a cyclic ordering of the edges coming into each vertex.  Then $(G, F)$ is obtained from $\Theta$ by blowing up each vertex into a circle, and distributing the edges incident to the vertex around the circle in their cyclic order.  The forest $F$ consists of  a maximal subtree of each of  the circles.    
\end{proof}

The cocycle $\mu_\Theta$ corresponds to a class 
 $$[\mu_\Theta]\in H_{2r-2+a}(\Out (F_{r+a});\Q)$$
 where $a$ is the number of vertices of $\Theta$ and $r$ is the rank of $\pi_1(\Theta)$. The homology $H_i(\Out(F_r);\Q)$ is independent of $r$ for large $r$ \cite{HaVo2}, but these classes are not in the  stable range.  
 
It remains an open question whether the cocycles $\mu_\Theta$ correspond to non-trivial homology classes.

\subsection{Cocycles parameterized by graphs with two types of trivalent vertices}
Let $c\g^{\prime\prime}$ be defined like $c\g^\prime$, except
that the trivalent vertices of a graph come in two types,  which we call   type A and type B. The graphs are oriented by cyclically ordering the edges at each type A vertex 
(If all vertices are type A, this is equivalent to the standard notion of orientation, see \cite[Proposition 2]{exposition}, since there are no even-valent vertices.) 
Now $\tr\co f\g^\prime\to c\g^\prime$ has an evident generalization to a map
$\Tr^{\g^\prime}\co f\g^\prime\to c\g^{\prime\prime}$ as follows. This map is non-trivial on $(G, F)$ only if each component of $F$ is either a linear tree $T$ with an even number of edges  whose endpoints are joined by an edge $e_T$ of $G-F$, or is a single vertex.  The map is defined as before by collapsing the cycles $T\cup e_T$  to oriented vertices (these are the type A vertices in the image).  The components of $F$ which are single vertices become type B vertices in the image. There is also a factor of $(-2)^a$ where $a$ is the number of type A vertices in the image.

If $\Theta\in c\g^{\prime\prime}$ has $a$ vertices of type A,  $b$ vertices of type B, and rank of $\pi_1(\Theta)=r$,
then composition of $\Tr^{\g^\prime}$ with projection of $c\g^{\prime\prime}
$ onto the subspace spanned by $\Theta$ gives a cocycle which corresponds to a class  
$$[\mu_\Theta]\in H_{2r-2+a-b}(\Out (F_{r+a});\Q).$$
Note that if there are enough type B vertices, these classes lie in the stable range.

It is natural to wonder whether the process of modifying an odd-valent graph by adding new edges
along type B vertices placed in the middles of existing edges is related to stabilization.

\section{Bordification}

In this section we describe another construction of cycles   in a different   chain complex which also computes the rational homology of $\Out(F_r)$.   This chain complex arises from the  Bestvina-Feighn bordification of Outer space \cite{BeFe},  and leads to a host of new  cycles, which also lift to rational cycles for $\Aut(F_r)$.  We conjecture that these new cycles include the Morita cycles and their generalizations.

\subsection{Bordified outer space}
Outer space is a contractible space  on which $\Out(F_r)$ acts with finite stabilizers \cite{CuVo}, so that
the quotient of Outer space by the group action is  a rational classifying space for $\Out(F_r)$.
Outer space is an equivariant union of open simplices, but some of the simplices have missing faces, so the  quotient is not compact and it is not possible to read off a chain complex for the rational homology of $\Out(F_r)$ from this simplicial structure.   One way to surmount this difficulty is to consider the so-called \emph{spine} \cite{CuVo},   which is a locally finite simplicial complex onto which Outer  space deformation retracts. A complementary method of obtaining a locally finite CW complex by attaching extra cells to the boundary was introduced in \cite{BeFe}. The result is called the {\it bordification}, and we denote it by $\B_r$.

Cells of $\B_r$ correspond to {\it marked filtered $r$--graphs}.   Define a \emph{core graph} to be a graph with vertices of valence at least $2$ and no separating edges. 
Then  a \emph{filtered $r$--graph} $\mathbf G$ is a connected graph $G$ of rank $r$, with all vertices of valence at least 3, together with a (possibly empty) chain $G_1\subsetneq \cdots \subsetneq G_{k-1}$ of proper core subgraphs of $G$. A {\it marking} of a filtered $r$--graph is a homotopy equivalence $\phi\co G\to R_r$, where $R_r$ is a fixed   wedge of $r$ circles (also called a rose with $r$ petals).
Two marked filtered $r$--graphs $( \phi, G_1\subset\cdots\subset G_k=G)$ and $(\phi',G_1'\subset\cdots\subset G_k'=G')$ are equivalent if there is a cellular isomorphism $f\co G\to G'$ preserving the filtrations, with $\phi'
\circ f \simeq\phi$.

There is one (homotopy) cell of $\B_r$ for each marked filtered $r$--graph. The dimension of the cell is equal to the number of edges in the graph minus the length of the filtration. The codimension $1$ faces of a cell are formed by either inserting a new core graph in the filtration, or by collapsing an edge; however, one is not allowed to collapse a loop (which would result in a graph with smaller rank fundamental group) or an edge which is equal to $G_i -G_{i-1}$ for some $i$ (which would decrease the length of the filtration).   

For more about $\B_r$, see \cite{BeFe, cut}. 

\subsection{Bordified auter space}

A space similar to Outer space was introduced in \cite{HaVo2} for $\Aut(F_r)$, and is sometimes referred to as ``Auter space."  The definition and auxiliary constructions are entirely analogous to those of Outer space, except that marked graphs have basepoints.  We will use $\B^\dagger_r$ to denote the bordification of Auter space.  

Define a \emph{basepointed filtered $r$--graph}, $\mathbf G$ to be a filtered $r$--graph which has a specified basepoint. All vertices of the whole graph are of valence at least three except the basepoint, which is allowed to be of valence $2$. A marking of a basepointed filtered $r$--graph is a homotopy equivalence to $R_r$ which preserves the basepoint. Two marked filtered graphs are equivalent if they differ by a basepoint preserving cellular isomorphism.
There is  one cell in $\B^\dagger_r$ for every basepointed filtered $r$--graph, and the cells fit together as in $\B_r$.

\subsection{Chain complexes}

Bestvina and Feighn prove that $\B_r$ is contractible, and that $\Out(F_r)$ acts properly discontinuously; the same is true for and $\B_r^\dagger$ and the action of $\Aut(F_r)$, so that:
\begin{thm}  {\rm \cite{BeFe}}

\noindent\begin{enumerate}
\item $H_i(\Out(F_r);\mathbb Q)\cong H_i(\B_r/\Out(F_r);\mathbb Q)$
\item $H_i(\Aut(F_r);\mathbb Q)\cong H_i(\B^\dagger_r/\Aut(F_r);\mathbb Q)$
\end{enumerate}
\end{thm}

In order to do homology computations, we need an explicit description of the generators and differentials of the chain complex for $\B_r/\Out(F_r)$.   We now describe a chain complex $\b_*$ and show that it can be used to compute the homology of  $\B_r/\Out(F_r)$.  (An analogous complex  $\b_*^\dagger$ computes the homology of  $\B_r^\dagger/\Aut(F_r)$).  The main issue is translating the notion of orientation of  a cell in the bordification to a notion of orientation on the associated filtered $r$--graph.

We define an {\it orientation} on a  filtered graph  $\mathbf G=\{G_1\subset\cdots\subset G_{k}=G\}$ to be an ordering of the edges of $G$, up to even permutation.  If $\mathbf G\langle C\rangle$ is obtained from $\mathbf G$ by inserting a core subgraph $C$ into the filtration, then $\mathbf G\langle C\rangle$ naturally inherits an orientation from $\mathbf G$.  If $\mathbf G_e$ is obtained from $\mathbf G$ by collapsing an edge $e$ of $G$, we orient $\mathbf G_e$ by simply leaving out $e$ in the ordering of the edges.  

Define the {\it degree}  of $\mathbf G$ to be equal to $e(G)-k$, where $e(G)$ is the number of edges of $G$.   Let $\b_*$ be the rational graded vector space spanned by isomorphism classes of oriented filtered $r$-graphs, modulo the relation that $(\mathbf{G},-or)=-(\mathbf{G},or)$, where the grading on  $\b_*$ is  by degree.  The differential $d\co \b_i\to\b_{i-1}$ is the sum of two differentials $d=d_E+d_F$. Here
$$d_F(\mathbf{G})=\sum_C(-1)^{f(C)} \mathbf G\langle C\rangle,$$  where the sum is over all core graphs $C$ which can be inserted into the filtration, and $f(C)$ is the position of  $C$ in the filtration after it is inserted, and 
$$d_E(\mathbf{G})=\sum_e (-1)^{n(e)+k-1}\mathbf G_e,$$ 
where the sum is over all edges $e$ which can be collapsed while preserving both the rank of $G$ and the number of graphs in the filtration, and $n(e)$ is the position of the edge $e$ in the ordering.

We remark that the chain complex $\b_*$ can be considered as a double complex,  bi-graded by the number of vertices in the whole graph and the number of graphs in the filtration. The differential $d_E$ decreases the number of vertices by 1, and the differential $d_F$ increases the number of graphs in the filtration by 1.  

An identical construction gives a bi-graded chain complex $\b^\dagger_*$ spanned by basepointed filtered $r$--graphs.  

\begin{prop}
\noindent\begin{enumerate}
\item $H_i(\b_*)\cong H_i(\B/\Out(F_r);\mathbb Q)$.
\item $H_i(\b^\dagger_*)\cong H_i(\B^\dagger/\Aut(F_r);\mathbb Q)$
\end{enumerate}
\end{prop}
\begin{proof}  

A cell in the bordification contributes a generator to the chain complex of the quotient $\B_r/\Out(F_r)$
if and only if all elements of its stabilizer preserve its orientation (see  \cite[Sec. 3]{HaVo} for details in a similar situation).   Thus to identify the generators of the chain complex, as well as to verify the signs in the definition of the boundary operators $d_E$ and $d_F$, we need to 
 to show how our notion of orientation of a graph is related to the standard notion of orientation of a cell in the bordification.
 
We first describe the orientation of a cell.  Given a filtered graph $\mathbf G=\{G_1\subset\cdots\subset G_{k}=G\}$ the open cell corresponding to it, denoted $\Sigma_{\mathbf G}$, breaks up as a Cartesian product $\Sigma_{G_1}\times \Sigma_{G_2/G_1}\times\cdots\times\Sigma_{G/G_{k-1}}$\cite[Proposition 2.15]{BeFe}. Each term in this  product  is parameterized by varying lengths on the edges of $G_i-G_{i-1}$ subject to the constraint that the sum of these lengths is $1$. Expanding the length of an edge $e$ of $G_i-G_{i-1}$ while uniformly shrinking the other edges gives rise to a vector $v_e$ in the tangent bundle of $\Sigma_{G_i/G_{i-1}}$ and hence of the entire cell $\Sigma_G$.   Because the total volume of $G_i/G_{i-1}$ is $1$, we have $\sum_{e\in G_i-G_{i-1}} v_e=0$.
An ordering $(e_1,\ldots,e_n)$ of the edges in $G_i-G_{i-1}$ gives an orientation $v_{e_2}\wedge\ldots
\wedge v_{e_n}$, and the reader can check that a permutation of $(e_1,\ldots,e_n)$ changes the induced orientation by the sign of the permutation.

We need to determine how this orientation behaves upon passing to a codimension 1 face. Let $\nu$ be an inward-pointing normal vector to such a face. By convention,  an orientation $\beta$ of the face is induced by the orientation $\alpha$ of the interior if $\nu\wedge\beta=\alpha$. 

There are two types of faces, the faces $\Ge$ coming from collapsing an edge $e$, and the faces $\GC$ coming from inserting a core graph $C$ into the filtration.   Each of $\Ge$ and $\GC$ has an orientation inherited from the orientation of $\mathbf G$.   We will show that the orientation on $\mathbf G_e$ induced from the interior is $(-1)^{n(e)+f(e)-1}$ times the orientation inherited from $\mathbf G$, where $n(e)$ is the position of $e$ in the ordering, and $f(e)$ is the stage of the filtration at which $e$ first appears.  For $\GC$,  we will show that the orientation induced from the interior is $(-1)^{f(C)+|C|+1}$ times the orientation inherited from $\mathbf G$,  where $f(C)$ is the position of $C$ in the filtration, and $|C|$ is the number of edges in $C$. In both of these sign conventions, changing the ordering of the edges of the graph $\mathbf G$
changes the orientation induced from the interior by the sign of the permutation. Thus, to establish these signs, it suffices to consider any fixed ordering of the edges in $\mathbf G$.

For notational convenience, then, assume that our ordering of the edges of $\mathbf G$ is such that the edges in $G_i$ lie before those in $G_j$ for $i<j$. This gives rise to an orientation
$\alpha_1\wedge\cdots\wedge\alpha_k$ where $\alpha_i=v_{e_2}\wedge\cdots\wedge v_{e_n}$ and $(e_1,\ldots,e_n)$ is the ordering of the edges in $G_i-G_{i-1}$.

Now let us examine what happens when we pass to a face by contracting an edge $e$.
In this case the normal vector to the face is $v_e$.
Let $f=f(e)$, ie, $e$ is in $G_f-G_{f-1}$.  Since we are allowed to contract $e$, $e$ can't be the only edge of $G_f-G_{f-1}$, so we may assume that $e$ occupies the second spot in the ordering. The orientation of the interior  of the cell is then
$$\alpha_1\wedge \alpha_2\wedge\cdots\wedge\alpha_{f-1}
\wedge v_e\wedge \alpha_f^\prime\wedge
 \alpha_{f+1}\wedge\cdots\wedge\alpha_k,$$ 
 where $\alpha_f^\prime$ is the wedge of vectors corresponding to the edges in $G_f-G_{f-1}$ coming after $v_e$. Bringing the normal vector forward, we get
 $$(-1)^{\epsilon}v_e\wedge\alpha_1\wedge \alpha_2\wedge\cdots\wedge\alpha_{f-1}\wedge
 \alpha_f^\prime\wedge
 \alpha_{f+1}\wedge\cdots\wedge\alpha_k,$$
where $\epsilon =|\alpha_1|+\cdots+|\alpha_{f-1}|$. We now observe that $\epsilon =(n(e)-2)-(f-1)=n(e)-f-1$, which is equal to  $n(e)+f(e)-1$ mod 2, as claimed.  

 This gives rise to a boundary operator $d^\prime_E$ defined by
$$d^\prime_E(\mathbf G) =\sum_{e\in \mathbf G} (-1)^{n(e)+f(e)-1}\mathbf G_e,$$where $n(e)$ is the number of $e$ according to the ordering, and $f(e)$ is the stage in the filtration that $e$ lies.

Next we examine what happens when we pass to a face by inserting a core graph into the filtration.
If $C$ is a core graph that can be inserted in the filtration of $\mathbf G$, define
$$\mathbf G\langle C\rangle= G_1\subset\cdots\subset G_{f-1}\subset C\subset G_f\subset
\cdots\subset G,$$
where $f=f(C)$.
Suppose that $G_f-G_{f-1}$ has ordered edges  $(e_{1},\ldots,e_N)$
and that $C-G_{f-1}$ has edges $(e_{1},\ldots,e_M)$.
In this case, the normal vector is gotten by expanding all the edges in $C$ uniformly:
$\nu= v_1+\cdots +v_M=-(v_{M+1}+\cdots + v_N)$.
Let $\alpha_C= v_2\wedge\cdots\wedge v_M$, and $\alpha_f^\prime=v_{M+2}\wedge\cdots\wedge v_N$.
Now the orientation on the interior of our cell is 
\begin{align*}
&\alpha_1\wedge\cdots\wedge\alpha_{f-1}\wedge
\alpha_C\wedge v_{M+1}\wedge\alpha^\prime_f\wedge \alpha_{f+1}\wedge\cdots\wedge \alpha_k\\
&=\alpha_1\wedge\cdots\wedge\alpha_{f-1}\wedge
\alpha_C\wedge (-\nu) \wedge\alpha^\prime_f\wedge \alpha_{f+1}\wedge\cdots\wedge \alpha_k.
\end{align*}
Pulling $\nu$ to the front, we get a coefficient of $(-1)^{\epsilon}$ where
$\epsilon = |\alpha_1|+\cdots+|\alpha_{f-1}|+|\alpha_C| + 1$, which is  
equal to  $|C|-f+1$.

This gives  a boundary operator $d^\prime_F$  given by 
$$d^\prime_F(\mathbf G)=\sum_{C} (-1)^{f(C)+|C|+1}\mathbf G \langle C\rangle.$$
Finally we introduce an automorphism of $\b_*$ which takes $d^\prime_E, d^\prime_F$ to the boundary operators $-d_E, -d_F$. 
Namely, if $\mathbf G=\{G_1\subset G_2\subset\cdots\subset G_k=G\}$ define $A (\mathbf G) = (-1)^{\sum_{i=1}^k|G_i|}\mathbf G$.  The same automorphism also works for $\b^\dagger_*$.
\end{proof}

\subsection{Constructing cycles from \ags}

In this section we will construct a map from  a set of combinatorial objects to $H_*(\b_*)$. Included in these objects will be the set of oriented odd-valent graphs, and we conjecture that the images of these are the same as the Morita cycles $[\mu_{\Theta}]$  we constructed earlier.

\eject
\begin{definition} 
\noindent\begin{itemize}
\item An  \emph{\ag}\! is a connected graph with vertices of valence $\geq 3$  which are of two types: A and B.
An \ag is called \emph{odd} if all of the type A vertices are of odd valence.
\item An \emph{orientation} of an odd \ag\! is determined by an ordering of the half-edges incident to each type A vertex. Reordering the half-edges will flip the orientation if the permutation is odd and fix it if the permutation is even.
\item $\m$ is the vector space spanned by oriented odd \ags modulo the relations
$(X,or)=-(X,-or)$.
\end{itemize}
\end{definition}

Given an oriented odd AB-graph $X$, we will construct a set of oriented filtered graphs, corresponding to oriented cells in the bordification $\B$.  A suitable linear combination of these cells will form a cycle, giving us a map
$$\Phi\co\m\to H_*(\b_*).$$
To define $\Phi(X)$, we first make some choices, and will show later that $\Phi(X)$ is independent of these choices.  Choose 
\begin{itemize}
\item a maximal tree $T$ of $X$,
\item an ordering $\o_1=(s_1,\ldots, s_a)$ of the type A vertices of $X$,
\item an ordering $\o_2=(e_1,\ldots, e_k)$ of the edges of $X-T$.
\end{itemize}

\begin{figure}
 \begin{center}
\includegraphics[width=2in]{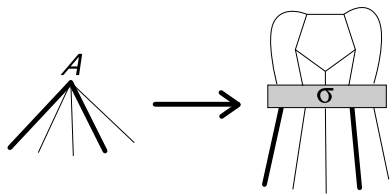}
\caption{Inserting a pentagon at a $5$--valent type A vertex. The boxed $\sigma$ indicates a bijection from the vertices of the pentagon to the edges adjacent to $v$. The heavy lines indicate the tree $T$.
}\label{pentagon}
\end{center}
\end{figure}

 To every type A vertex $s$ of $X$, associate a $v$--gon  $\Delta_s$, where $v$ is the valence of $s$, with the edges of $\Delta_s$ ordered cyclically $1,\ldots, v$.   If we fix bijections $\iota_s$ from the edges coming into $s$ to the vertices of  $\Delta_s$, we can construct a  new graph by replacing each $s$  by the associated $v$--gon $\Delta_s$, glued according to $\iota_s$.  Figure~\ref{pentagon} depicts this construction at a $5$--valent vertex.
  We then collapse the image of  $T$ (which is no longer a maximal tree) to obtain a graph $Y=Y[\iota_{s_1},\ldots,\iota_{s_a}]$.  To orient $Y$, we use the order $\o_1$ to totally order the edges in the union of the (edge-ordered) polygons $\Delta_s$, then order the remaining edges according to the order $\o_2$.  We filter $Y$ in an analogous manner, by successively adding the polygons, then the edges of $X-T$, in the order given by $\o_1$ and $\o_2$:
 \begin{align*}
\Delta_{s_1}&\subset\Delta_{s_1}\cup\Delta_{s_2}\subset\cdots\subset\Delta_{s_1}\cup\cdots\cup \Delta_{s_a}\subset
\Delta_{s_1}\cup\cdots\cup \Delta_{s_a}\cup e_1\\
&\subset\Delta_{s_1}\cup\cdots\cup \Delta_{s_a}\cup e_1\cup e_2\subset\cdots\subset \Delta_{s_1}\cup\cdots\cup \Delta_{s_n}\cup e_1\cup \cdots \cup e_k=Y
\end{align*}
The ordering on the edges of each $\Delta_s$ naturally gives an ordering $1,\ldots,v$ of the vertices.  The edges of $X$ coming into $s$ can also be ordered $1,\ldots, v$ using the orientation data for $X$.  Thus the set $I_s$ of all bijections from the  incoming edges at $s$ to the vertices of $\Delta_s$ is naturally identified with the set of permutations of $v$ letters.   For $\iota\in I_s$ let $\epsilon(\iota)$ be the sign of the corresponding permutation, and define
$$\Phi(X)=\sum_{\iota_j\in I_{s_j}} \epsilon(\iota_1)\cdots\epsilon(\iota_n) Y[\iota_1\cdots \iota_n].$$
Note that it is critical here that all type A vertices have odd valence. The dihedral group $D_v$ acts on $I_s$, and the isomorphism type of the graphs $Y$  is constant on each orbit. When $v$ is even the corresponding terms of $\Phi(X)$ all cancel, giving $\Phi(X)=0$,  whereas if $v$ is odd, they all have the same sign.

\begin{proposition}
$\Phi(X)$  is a cycle.
\end{proposition}
\begin{proof}
The filtration on each $Y$  is maximal,  so the boundary operator $d_F$ vanishes on $\Phi(X)$. The  second boundary operator $d_E$ also vanishes. To see why first notice that the only edges of $Y$ available for collapsing are in the polygons $\Delta_s$.  When we contract an edge $e$ in $\Delta_s$, the resulting filtered graph $Y_e$ is the same as we get by replacing the bijection $\iota_s$ by a bijection which transposes the two edges incident to $e$.  Thus the terms cancel in pairs.
\end{proof}

\label{independent}
\begin{proposition}  $\Phi(X)$ is well-defined, ie, $\Phi(X)$  is independent of the choice of maximal tree $T$, the ordering $o_1$ of type A vertices and the ordering $o_2$ of edges of $X-T$.
\end{proposition}
\begin{proof}
First we address varying the tree.
Suppose one replaces an edge $e$ of $T$ with an edge $e^\prime$ of $X-T$ to obtain a new maximal tree $T^\prime$.  Give $e'$ the same place in the ordering of edges of $X-T'$ as $e$ occupied in the ordering of $X-T$,  and form a new cycle $\Phi'(X)$.  Now consider the sum of oriented filtered graphs formed by inserting polygons at each type A vertex of $X$ using all possible bijections as above, and collapsing all the edges in $T$  {\it other} than $e$. The orientation and filtration of the summands are defined as  in the definition of $\Phi(X)$, by first including the polygons one at a time and then the remaining edges, except that $e$ and $e^\prime$ are added at the same time in the filtration, and $e$ comes before $e'$ in the orientation. This defines a new chain, whose boundary is the difference of  $\Phi(X)$ and $\Phi'(X)$, so they are homologically equal. (The minus sign comes from the fact that we have to switch the order of $e$ and $e^\prime$.)
The argument is finished by the observation that one can get between any two maximal trees by replacing one edge at a time.

Now we vary $\o_1$. Suppose $\o_1$ and $\o^\prime_1$ differ by a consecutive transposition (for simplicity, say $1\leftrightarrow 2$), and that $\Phi'(X)$ is constructed using $\o^\prime_1$. In each term of $\Phi(X)$,  $\Delta_{s_1}\cup\Delta_{s_2}$ is either a disjoint union of two circles or it is a wedge of two circles. In either case the $\Delta_{s_i}$ are the only proper core subgraphs of their union. 
Consider the chain $z$ which is the sum of filtered graphs  $\Phi(X)$, but with the first term of the filtration truncated in each summand:
$$ \Delta_{s_1}\cup\Delta_{s_2}\subset \cdots.$$
The boundary $d_E(z)$  vanishes, whereas $d_F(z)$ either inserts $\Delta_{s_1}$ or $\Delta_{s_2}$ first, the latter with the sign coming from moving the edges of $\Delta_{s_1}$ past those of $\Delta_{s_2}$. Since both have an odd number of edges, this sign is $-1$. Thus, as before, we see that the difference $\Phi(X)-\Phi^\prime(X)$ is a boundary.

The argument for replacing $\o_2$  by $\o^\prime_2$ is very similar.
\end{proof}

If  $X$ has  $a$ vertices of type A, $b$ vertices of type $B$ and fundamental group of rank $n$, and if $V_B$ denotes the sum of the valences of the type B vertices, then 
$$\Phi(X)\in H_{2n-2+a+2b-V_B}(\Out (F_{n+a});\Q).$$
We remark that when all of the type B vertices are of valence $3$, then the degree is the same as that of $[\mu_X]$.

If all type $A$ vertices of $X$ happen to be trivalent, then all the ways of inserting the triangles at each vertex are equivalent, so that $\Phi(X)$ can be represented as a single filtered graph, corresponding to a single cell of $\B$.
For example, the image of the theta graph $\Theta_3$ is a multiple of the class pictured in Figure~\ref{theta}.

\begin{figure}
\begin{center}
\includegraphics[width=3in]{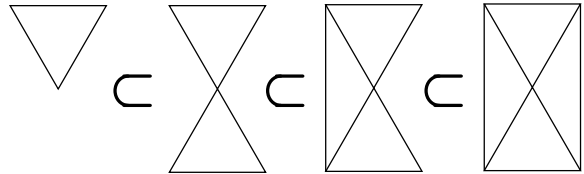}
 \caption{$\Phi(\Theta_3)$ in the bordified chain complex.}\lbl{theta}
 \end{center}
 \end{figure}

A straightforward but messy calculation shows that $\Phi(\Theta_3)$ is non-trivial, so it generates $H_4(\Out(F_4);\Q)$ and therefore must be a multiple of the Morita cycle $[\mu_{\Theta_3}]$, we constructed earlier.  This supports the following conjecture:

\begin{conjecture}
If $X$ is an odd-valent graph, then $[\mu_X]$ is a multiple of $\Phi(X)$. 
\end{conjecture}

We conclude by noting that the set $\{\Phi(X)\}$ of homology classes is definitely not independent. For example, if $X$ is a graph with no type A vertices and only type B vertices then $\Phi(X)$ always represents a nonzero class in $H_0(\Out(F_r);\Q)$ where $r$ is the rank of the  fundamental group of $X$. However it seems reasonable to conjecture that the set of graphs with only type A vertices gives rise to a (largely) independent set.

\subsection{Classes for ${\rm Aut}(F_r)$}
Let $\m^\dagger$ be the space of odd oriented \ags with a specified basepoint. 
Then there is a map
$$\Phi^\dagger\co \m^\dagger\to H_*(\b_*^\dagger),$$
defined similarly to the map $\Phi$ of the previous section.
The basepoint of the \ag descends to a basepoint of the image.
This is straightforward if it is a type B basepoint, but it also works for type A baspeoints by choosing the basepoint anywhere on the polygon to be inserted in the vertex. The choice is irrelevant since we are averaging over all ways to insert it.

For degree reasons, the nontrivial class in $H_7(\Aut(F_5);\Q)$ discovered by Gerlits is not in the image of $\Phi^\dagger$.

\end{document}